\documentclass[11pt]{amsart}
\usepackage[english]{babel}
\usepackage[dvips]{graphicx}
\usepackage{amsmath,amsfonts,amsthm,amssymb,dsfont}
\usepackage[usenames, dvipsnames]{color}

\newtheorem{theo}{Theorem}[section]

\newtheorem{lem}{Lemma}[section]

\newtheorem{rem}{Remark}[section]
\newtheorem{Example}{Example}[section]

\newcommand{\secondmomentpideltaD}{A1}
\newcommand{\constantbalanceH}{A2}
\newcommand{\constantbalanceD}{A'2}
\newcommand{\irreducibleH}{A3}

\title{Interacting generalized Friedman's urn systems}

\author[G. Aletti]{Giacomo Aletti}
\address[G. Aletti]{ADAMSS Center, Universit\`a degli Studi di Milano, 20133 Milano, Italy}
\email{giacomo.aletti@unimi.it}

\author[A. Ghiglietti]{Andrea Ghiglietti}
\address[A. Ghiglietti]{Universit\`a degli Studi di Milano, 20133 Milano, Italy}
\email{andrea.ghiglietti@unimi.it}

\date{\today}

\begin{document}

\maketitle

\begin{abstract}
We consider systems of interacting \textit{Generalized Friedman's Urns} (GFUs) having irreducible mean replacement matrices.
The interaction is modeled through the probability to sample the colors from each urn,
that is defined as convex combination of the urn proportions in the system.
From the weights of these combinations we individuate subsystems of urns evolving with different behaviors.
We provide a complete description of the asymptotic properties of urn proportions in each subsystem
by establishing limiting proportions, convergence rates and Central Limit Theorems.
The main proofs are based on a detailed eigenanalysis and stochastic approximation techniques.
\end{abstract}



\section{Introduction}\label{introduction}

The stochastic evolution of systems composed by elements
which interact among each other has always been of great interest in several areas of application,
e.g. in medicine a tumor growth is the evolution of a system of interacting cells~\cite{YakTso},
in socio-economics and life sciences a collective phenomenon reflects the result of the interactions among the individuals~\cite{NalParTos},
in physics the concentration of certain molecules within cells varies over time due to interactions between different cells~\cite{Pri}.
In the last decade several models have been proposed in which the elements of the system are represented by urns containing balls of different colors,
in which the urn proportions reflect the status of the elements,
and the evolution of the system is established by studying the dynamics at discrete times of
this collection of dependent urn processes.
The main reason of this popularity is concerned with the urn dynamics, which is
(i) suitable to describe random phenomena in different scientific fields (see e.g.~\cite{JohKot}),
(ii) flexible to cover a wide range of possible asymptotic behaviors,
(iii) intuitive and easy to be implemented in several fields of application.

The dynamics of a single urn typically consists in a sequential repetition of a \emph{sampling phase}, when a ball is sampled from the urn,
and a \emph{replacement phase}, when a certain quantity of balls is replaced in the urn.
The basic model is the P\'{o}lya's urn proposed in~\cite{EffPol}:
from an urn containing balls of two colors, balls are sequentially sampled and then replaced in the urn with a new ball of the same color.
This updating scheme is then iterated generating a sequence of urn proportions whose almost sure limit is random and Beta distributed.
Starting from this simple model, several interesting variations have been suggested by considering different distributions in
the sampling phase, e.g.~\cite{HilLanSud1,HilLanSud2}, or in the replacement phase, e.g.~\cite{Ath,Fri,Pem}.
In a general $K$-colors urn model,
the color sampled at time $n$ is usually represented by a vector $\mathbf{X}_n$
such that $X_{i,n}=1$ when the sampled color is $i\in\{1,{\ldots},K\}$, $X_{i,n}=0$ otherwise;
the quantities of balls replaced in the urn at time $n$ are typically
defined by a matrix $D_n$ such that $D_{ki,n}$ indicates the number of balls of color $k$ replaced in the urn when the color $i$ is sampled.
Considering $\{D_n;n\geq1\}$ as an i.i.d. sequence,
a crucial element to characterize the asymptotic behavior of the urn
is the mean replacement matrix $H:=\mathbf{E}[D_n]$, typically called \emph{generating matrix}.

The class of urn models considered in this paper is commonly denoted by
\emph{Generalized Friedman's Urn} (GFU).
The GFU model was introduced in~\cite{Fri} and its extensions and their asymptotic behavior have been studied in
several works, see e.g.~\cite{AthKar2,BaiHu,BaiHu2,Smy}.
The GFU considered in this paper is characterized by a non-negative irreducible generating matrix $H$
with average constant balance, i.e. the columns of $H$ sum up at the same constant, $\sum_{k=1}^KH_{ki}=c>0$ for any $i\in\{1,{\ldots},K\}$,
which implies that its maximum eigenvalues $\lambda_{\max}(H)=c$ has multiplicity one.
The irreducibility of $H$ distinguishes the GFU from the \emph{Randomly Reinforced Urn} (RRU) model,
which includes the classical P\`{o}lya's Urn, whose replacement matrix is diagonal: when the color $i$ is sampled,
the GFU replaces in the urn more colors
following the distribution of the $i^{th}$ column of $D_n$ while the RRU only adds balls of colors $i$;
hence, the probability to sample color $i$ at next step is reinforced in the RRU,
while it may increase or decrease according to the current urn composition in the GFU.
The asymptotic behavior is in general very different:
in a GFU the urn proportion converges to a deterministic equilibrium identified by $H$ (see e.g.~\cite{AthKar2,BaiHu,BaiHu2,Smy}),
while in a RRU the limit is random and its distribution depends on the initial composition (see e.g.~\cite{AleMaySec1,AleMaySec2,DurFloLi}).

The model proposed in this paper is a collection of $N\geq1$ GFUs that interact among each other during the sampling phase:
the probability to sample a color $i$ from an urn $j$ is a convex combination of the urn proportions of the entire system.
Hence, a crucial role to describe the system dynamics is played by the \emph{interacting matrix} $W$ made by the weights of those combinations.
Since the asymptotic properties of the single GFUs are typically determined by the corresponding generating matrices $\{H^j;1\leq j\leq N\}$
and since the interaction among them is ruled by $W$,
the system dynamics has been studied by defining a new object $\mathbf{Q}$ that merges the information contained in $\{H^j;1\leq j\leq N\}$ and $W$.
From the analysis of the eigen-structure of $\mathbf{Q}$,
we are able to establish the convergence and the second-order asymptotic behavior of the urn proportions in the entire system.
Hence, this paper extends the theory on GFU models in the sense that,
in the special case of no interaction, i.e. $W=I$,
the results presented for the system reduce to the well-known results for a single GFU.

Several interacting urn models have been proposed in the last decade, especially for RRU models.
An early work is represented by~\cite{PagSec} that considered a collection of two-colors RRU in which
the replacements depend on the colors sampled in the rest of the system and hence
the sequence $\{D_n;n\geq1\}$ is not i.i.d.
Consequently, the interaction in~\cite{PagSec}
is modeled through the definition of $D_n$, instead of $\mathbf{X}_n$ as in our model.
A completly different updating rule has been used in the two-color urn model proposed in~\cite{MarVal},
in which sampling color 1 in the urn $j$ increases the composition of color 1 in the urn $j$,
while sampling color 2 increases the composition of color 2 in the neighbor urns $i\neq j$ and the urn $j$ comes back to the initial composition.
Asymptotic properties for this system have been obtained in~\cite{MarVal}
where there is no convergence of the urn proportions.
Other models in which the interaction enters in the replacement matrices are for instance~\cite{BenBenCheLim,CheLuc,CirGalHus}.

Recently there have been more works concerning urn systems in which the interaction is modeled
through the sampling probabilities as in our model.
They differ from this paper since all of them consider RRUs and the interaction is only modeled as mean-field interaction tuned by a parameter $\alpha\in(0,1)$,
i.e. the urns interact among each other only through the average composition in the entire system.
As a consequence, their asymptotic results lead to the \textit{synchronization} property in which all the urn proportions of the system
converge to the same random limit.
In particular, in~\cite{Lau,LauLim} the asymptotic behavior of the urn system has been studied
for a model that defines the sampling probabilities through the exponential of the urn compositions.
In~\cite{CriDaiPraMin,DaiPraLouisMin} the sampling probabilities are defined directly using the urn compositions
and, in addition, the synchronization property has been proved;
moreover, different convergence rates and second-order asymptotic distributions for the urn proportion have been established
for different values of the tuning parameter $\alpha$.
Since we consider GFU models the asymptotic results established in this paper are totally different
from those proved in~\cite{CriDaiPraMin,DaiPraLouisMin}, e.g.
our limiting proportions are not random and they do not depend on the initial compositions.

It is also significant to highlight that this work allows a general structure for the urn interaction,
which reduces to the mean-field interaction only for a particular choice of the interacting matrix $W$.
Moreover, from the analysis of the structure of $W$ we are able to individuate subsystems of urns evolving with different behaviors
(see Subsection~\ref{subsection_W interaction}):
(i) the \textit{leading systems}, whose dynamics are independent of the rest of the system
and (ii) the \textit{following systems}, whose dynamics ``follow'' the evolution of other urns of the system;
in the special case of irreducible interacting matrix, which includes the mean-field
interaction considered in~\cite{CriDaiPraMin,DaiPraLouisMin},
there is a unique leading system and no following systems.
These two classes of systems have been studied separately
(leaders in Section~\ref{section_leading} and followers in Section~\ref{section_follower}),
in order to provide an exhaustive description of the asymptotic behavior in any part of the system.
In fact, since different systems may converge at different rates,
a unique central limit theorem would not be able
to characterize the convergence of any urn proportion.
Hence, through a careful analysis on the eigen-structure of $\mathbf{Q}$ realized in Subsection~\ref{subection_in_out_follower},
we individuate the components of the urn processes in the system that actually ``lead'' or influence the following systems,
so that we can establish the right convergence rate and a non-degenerate asymptotic distribution for any subsystem.

A pivotal technique in the proofs consists in revisiting the dynamics of the urn proportions of the system in the stochastic approximation (SA) framework,
as suggested for the composition of a single GFU in~\cite{LarPag}.
To this end, the dynamics of the urn compositions of the same subsystems
have been reformulated into a recursive stochastic algorithm (see Section~\ref{section_stochastic_approximation}).
Then, the dynamics of the urn proportions have been properly modified to embed the processes
of the urn proportions into the whole suitable space $\mathbb{R}^K$
(see Subsection~\ref{subection_redefining_simplex_leading} and~\ref{subection_redefining_simplex_follower}).\\

The main results of the paper starts at Section~\ref{section_leading}.
The first part of the paper is a necessary formulation of the problem in its general form,
together with all the assumptions and notations that may appear tough at a first reading.
We provide a guiding Example~\ref{example_1}, that is recovered in the Example~\ref{example_2} and Example~\ref{example_3},
to help the reader to appreciate the main results, although not in all their depth.

More precisely, the structure of the paper is the following.
In Section~\ref{section_model_assumptions} we present model and main assumptions concerning the interacting GFU system.
Specifically, in Subsection~\ref{subection_model} we describe how the composition of the colors in each urn of the system evolves at any time $n\geq1$.
Then, in Subsection~\ref{subection_main_assumptions} the main assumptions required to establish the results of the paper are presented.
Subsection~\ref{subection_preliminary} contains a preliminary result.
Subsection~\ref{subsection_W interaction} is dedicated to analyze the structure of the interacting matrix and
hence to define the leading and the following subsystems that compose the entire system.

Section~\ref{section_stochastic_approximation} is concerned with the dynamics of the interacting GFU system
expressed in the stochastic approximation framework.
In particular, in Subsection~\ref{subection_notation_SA} we introduce the notation
that combines the composition of the urns in the same subsystem.
Then, in Subsection~\ref{subection_dynamics_SA} the dynamics of the urn proportions
in any subsystem is reformulated into a recursive stochastic algorithm.

Section~\ref{section_leading} and Section~\ref{section_follower} contain the main results of the paper.
In particular, Section~\ref{section_leading} is concerned with the asymptotic behavior of the leading systems:
the convergence of the urn proportions is established in Subsection~\ref{subection_first_order_leading}
and the corresponding CLT is presented in Subsection~\ref{subection_second_order_leading}.
Then, Section~\ref{section_follower} is focused on the asymptotic behavior of the following systems:
in Subsection~\ref{subection_first_order_follower} we present the result on the convergence of the urn proportions,
while in Subsection~\ref{subection_second_order_follower} we establish the relative CLT.

Section~\ref{section_further_extensions} contains a brief discussion on further possible extensions of the interacting GFU model.
The proofs of all the results presented in the paper are contained in Section~\ref{section_proofs}.
Finally, in Appendix we report basic results of stochastic approximation that have been used in the main proofs.

\section{Model Setting and main Assumptions}\label{section_model_assumptions}

Consider a collection of $N\geq1$ urns containing balls of $K\geq1$ different colors.
At any time $n\geq0$ and for any urn $j\in\{1,{\ldots},N\}$,
let $Y^j_{k,n}>0$ be the real number denoting the amount balls of color $k\in\{1,{\ldots},K\}$,
let $T^j_n:=\sum_{k=1}^KY^j_{k,n}$ be the total number of balls
and let $Z^j_{k,n}:=Y^j_{k,n}/T^j_n$ be the proportion of color $k$.

\subsection{Model}\label{subection_model}

We now describe precisely how the system evolves at any time $n\geq1$.
Denote by $\mathcal{F}_{n-1}$ the $\sigma$-algebra generated by the urn compositions of the entire system up to time $(n-1)$, i.e.
\begin{equation*}
\mathcal{F}_{n-1}\ :=\ \sigma\left(\ X^j_{k,t}, Y^j_{k,t},\ 1\leq j\leq N,\ 1\leq k\leq K,\ 1\leq t\leq n-1\ \right).
\end{equation*}
The dynamics of the system is described by two main phases: \emph{sampling} and \emph{replacement}.\\

\emph{Sampling phase}:
for each urn $j\in\{1,{\ldots},N\}$, a ball is virtually sampled and its color is represented as follows:
$X^{j}_{i,n}=1$ indicates that the sampled ball is of color $i$, $X^{j}_{i,n}=0$ otherwise.
We denote by $\tilde{Z}^{j}_{i,n-1}$ the probability to sample a ball of color $i$ in the urn $j$ at time $n$, i.e.
\begin{equation*}
\tilde{Z}^{j}_{i,n-1} :=\mathbf{E}\left[\ X^{j}_{i,n}\ |\ \mathcal{F}_{n-1}\ \right].
\end{equation*}
Given the sampling probabilities $\{\tilde{Z}^{j}_{i,n-1},1\leq j\leq N,1\leq i\leq K\}$,
the colors are sampled independently in all the urns of the system and hence,
for any $i\in\{1,{\ldots},K\}$, $X^{1}_{i,n},{\ldots},X^{N}_{i,n}$ are independent conditionally on $\mathcal{F}_{n-1}$.
We define the sampling probabilities as convex combinations of the urn proportions of the system.
Formally, for any urn $j\in\{1,{\ldots},N\}$ we introduce the weights $\{w_{jh};1\leq h\leq N\}$ such that $0\leq w_{jh}\leq1$ and
$\sum_{h=1}^Nw_{jh}=1$.
Thus, the probability to sample the color $i$ in the urn $j$ is defined as follows
\begin{equation}\label{def:Z_tilde}
\tilde{Z}^{j}_{i,n-1}\ :=\ \sum_{h=1}^Nw_{jh}Z^{h}_{i,n-1}.
\end{equation}\\

\emph{Replacement phase}:
after that a ball of color $i$ has been sampled from the urn $j$,
we replace $D^j_{ki,n}$ balls of color $k\in\{1,{\ldots},K\}$ in the urn $j$.
For any urn $j$ we assume that $\{D_{n}^j;n\geq1\}$ is a sequence of i.i.d. non-negative random matrices, where $D_{n}^j:=[D^j_{ki,n}]_{ki}$.
We will refer to $D^j_n$ as \emph{replacement matrix} and to $H^j:=\mathbf{E}[D_{n}^j]$ as \emph{generating matrix}.
Notice that $H^j$ are time-independent since $\{D_{n}^j;n\geq1\}$ are identically distributed
(see Subsection~\ref{section_further_extensions} for possible extensions).
Moreover, we assume that at any time $n$ the replacement matrix for the urn $j$, i.e. $D_{n}^j$,
is independent of the sampled colors, i.e. $\{X_{i,n}^j;1\leq j\leq N, 1\leq i\leq K\}$, and
independent of the replacement matrices of the other urns of the system, i.e. $D_{n}^{j_0}$ with $j_0\neq j$.\\

In conclusion, the composition of the color $k\in\{1,{\ldots},K\}$ in the urn $j\in\{1,{\ldots},N\}$ evolves at time $n\geq1$ as follows:
\begin{equation}\label{eq:dynamics_Y}
Y^j_{k,n}\ =\ Y^j_{k,n-1}\ +\ \sum_{i=1}^K D^j_{ki,n} X^j_{i,n}.
\end{equation}

\subsection{Main Assumptions}\label{subection_main_assumptions}

We now present the main conditions required to establish the results of the paper.
The first assumption is concerned with bounds for the moments of the replacement distributions.
Specifically, we require the following condition:
\begin{itemize}
\item[(\secondmomentpideltaD)] there exists $\delta>0$ and a constant $0<C_{\delta}<\infty$ such that, for any $j\in\{1,{\ldots},N\}$ and any $k,i\in\{1,{\ldots},K\}$,
$\mathbf{E}[(D^j_{ki,n})^{2+\delta}]<C_{\delta}$.
\end{itemize}
Note that $C_{\delta}$ does not depend on $n$ since $\{D_{n}^j;n\geq1\}$ are identically distributed.

The second assumption is the average constant balance of the urns in the system and
it is imposed by the following condition on the generating matrices $H^1,{\ldots},H^N$:
\begin{itemize}
\item[(\constantbalanceH)] for any $j\in\{1,{\ldots},N\}$ and $i\in\{1,{\ldots},K\}$, there exists a constant $0<c^j<\infty$ such that
$\sum_{k=1}^K H^j_{ki}\ =\  c^j$.
\end{itemize}
Note that (\constantbalanceH) guarantees that the average number of balls replaced in any urn is constant,
regardless its composition.
Assumption (\constantbalanceH) is essential to obtain the asymptotic configuration of the system, i.e. the limiting urn proportions.
The second-order asymptotic properties of the interacting urn system,  namely the rate of convergence and the limiting distributions,
are obtained by assuming a stricter assumption than (\constantbalanceH).
This condition is expressed as follows:
\begin{itemize}
\item[(\constantbalanceD)] for any $j\in\{1,{\ldots},N\}$, $i\in\{1,{\ldots},K\}$, $\mathbf{P}\left(\ \sum_{k=1}^K D^j_{ki,n}\ =\ c^j\ \right)=1$,
i.e. each urn is updated with a constant total amount of balls.
\end{itemize}
Naturally, (\constantbalanceD) implies (\secondmomentpideltaD) with $C_{\delta}=(\max_{j}\{c^j\})^{2+\delta}$.

Notice that, by defining $\widehat{Y}^j_{k,n}=(c^j)^{-1}Y^j_{k,n}$ and $\widehat{D}^j_{ki,n}=(c^j)^{-1}D^j_{ki,n}$ for all $n\geq1$, the urn dynamics in~\eqref{eq:dynamics_Y} can be expressed in the following equivalent form:
$$\widehat{Y}^j_{k,n}\ =\ \widehat{Y}^j_{k,n-1}+\sum_{i=1}^K\widehat{D}^j_{ki,n} \cdot X^j_{i,n},\ \ \ \ \ \ \
\widehat{Z}^j_{k,n-1}\ =\ \frac{\widehat{Y}^j_{k,n-1}}{\sum_{k=1}^K\widehat{Y}^j_{k,n-1}}\ =\
\frac{Y^j_{k,n-1}}{\sum_{k=1}^KY^j_{k,n-1}}\ =\ Z^j_{k,n-1}.$$
Therefore, from now on we will denote by $Y^j_{k,n}$ and $D^j_{ki,n}$ the normalized quantities $\widehat{Y}^j_{k,n}$ and $\widehat{D}^j_{ki,n}$
and hence (\constantbalanceH) and (\constantbalanceD) are replaced by the following conditions:
\begin{itemize}
\item[(\constantbalanceH)] for any $j\in\{1,{\ldots},N\}$ and $i\in\{1,{\ldots},K\}$, $\sum_{k=1}^K H^j_{ki}\ =\  1$.
\end{itemize}
\begin{itemize}
\item[(\constantbalanceD)] for any $j\in\{1,{\ldots},N\}$ and $i\in\{1,{\ldots},K\}$, $\mathbf{P}\left(\ \sum_{k=1}^K D^j_{ki,n}\ =\ 1\ \right)=1$.
\end{itemize}
In this case, (\constantbalanceD) implies (\secondmomentpideltaD) with $C_{\delta}=1$.

Finally, we consider \textit{Generalized Friedman's Urns} (GFUs) with irreducible generating matrices,
as expressed in the following condition:
\begin{itemize}
\item[(\irreducibleH)] for any $j\in\{1,{\ldots},N\}$, $H^j$ is irreducible.
\end{itemize}
This assumption will guarantee deterministic asymptotic configurations for the urn proportions in the system.

\begin{rem}
Extensions to non-homogeneous generating matrices $\{H_n;n\geq0\}$
are possible, as discussed in Section~\ref{section_further_extensions}.
In that case, assumption (\constantbalanceH) should be referred to the
limiting matrix $H^j:=a.s.-\lim_{n\rightarrow\infty}H^j_n$.
\end{rem}

\subsection{A preliminary result}\label{subection_preliminary}

The assumptions (\constantbalanceH) and (\constantbalanceD) on the constant balance are essential to obtain the following result on
the total number of balls in the urns of the system:
\begin{theo}\label{thm:total_number}
Let $T^j_n=\sum_{k=1}^KY^j_{k,n}$ be the total number of balls contained in the urn $j$ at time $n$.
Then, under assumptions (\secondmomentpideltaD) and (\constantbalanceH), $\{T^j_{n}-n;n\geq1\}$ is an $L^2$ martingale and,
for any $\alpha<1/2$,
\begin{equation}\label{eq:total_number_convergence}
n^{\alpha}\left(\frac{T^j_{n}}{n}-1\right)\ \stackrel{a.s./L^2}{\longrightarrow}\ 0.
\end{equation}
Moreover, under assumption (\constantbalanceD), $T^j_{n}=T^j_{0}+n$ a.s. and hence~\eqref{eq:total_number_convergence} holds for any $\alpha<1$.
\end{theo}

\subsection{The interacting matrix}\label{subsection_W interaction}

The interaction among the urns of the system is modeled through the sampling probabilities $\tilde{Z}^j_{i,n-1}$,
that are defined in~\eqref{def:Z_tilde} as convex combinations of the urn proportions of the system.
Formally, we denote by $W$ the $N\times N$ matrix composed by the weights $\{w_{jh},1\leq j,h\leq N\}$ of such linear combinations
and we refer to it as \emph{interacting matrix}.
We now consider a particular decomposition of $W$ that individuates subsystems of urns evolving with different behaviors.
The same decomposition is typically applied to the transition matrices in the context of discrete time-homogeneous Markov chains (see~\cite{Nor})
to characterize the state space.
For this reason, we first present the decomposition of $W$ in this 
framework,
and then we identify the subsystems of urns as the communicating classes of the state space.

Consider a discrete time-homogeneous Markov chain with state space $\{1,{\ldots},N\}$ and transition matrix $W$, i.e.
the element $w_{jh}$ now represents the probability of a Markov chain to move from state $j$ to state $h$ in one step.
It is well-known (see~\cite{Nor}) that the \emph{communication} relationship ($i\sim j$ if there exist $m,n\geq 0$ such that $[W^m]_{ij}>0$ and
$[W^n]_{ji}>0$) induces a partition of the state space into communicating classes (some of them are necessarily closed and recurrent,
with possibly some transient classes).
The maximum eigenvalue is $\lambda=1$ and its multiplicity reflects the number of recurrent classes.
Accordingly, let us denote by $\mathcal{L}$ the set of labels that identify the communicating classes, $n_L\geq1$ the multiplicity of $\lambda_{\max}(W)=1$,
and define the integers $n_F\geq 0$ and $1\leq r^{L_1}<\ldots<r^{L_{n_L}}<r^{F_1}<\ldots<r^{F_{n_F}}= N$
such that $W$ can be decomposed as follows (see \cite[Example 1.2.2]{Nor} for the analogous upper triangular case):
\begin{equation}\label{def:W}
\begin{gathered}
W:=\begin{pmatrix}
W^{L}&0\\
W^{FL}&W^{F}
\end{pmatrix},
\\
\begin{matrix}
W^L  :=
\left(
\begin{smallmatrix}
W^{L_1}&0&\ldots&0\\
0&W^{L_2}&\ldots&\ldots\\
\ldots&\ldots&\ldots&0\\
0&0&\ldots&W^{L_{n_L}}
\end{smallmatrix} \right)
\\
W^{FL}:=
\left(
\begin{smallmatrix}
W^{F_1L_1}&\ldots&W^{F_1L_{n_L}}\\
\ldots&\ldots&\ldots\\
W^{F_{n_F}L_1}&\ldots&W^{F_{n_F}L_{n_L}}
\end{smallmatrix} \right)
,&
W^F := \left(
\begin{smallmatrix}
W^{F_1}&0&\ldots&0\\
W^{F_2F_1}&W^{F_2}&\ldots&0\\
\ldots&\ldots&\ldots&\ldots\\
W^{F_{n_F}F_1}&W^{F_{n_F}F_2}&\ldots&W^{F_{n_F}}
\end{smallmatrix} \right).
\end{matrix}
\end{gathered}
\end{equation}
where:
\begin{itemize}
\item[(1)] for any $l\in\mathcal{L}$, $W^{l}$ is an $s^{l}\times s^{l}$ irreducible matrix,
where we let $s^{l}:=r^{l}-r^{l^{-}}$ and $l^{-}$ indicates the element in $\mathcal{L}$ that precedes $l$
(by convention $L_1^{-}\equiv \emptyset$ and $F_1^{-}\equiv L_{n_L}$);
\item[(2)] $\mathcal{L}:=\mathcal{L}_L\cup\mathcal{L}_F$, $\mathcal{L}_L:=\{L^1,{\ldots},L^{n_L}\}$ and $\mathcal{L}_F:=\{F^1,{\ldots},F^{n_F}\}$ are sets of labels that identify, respectively, recurrent and transient communicating classes in the state space ($\mathcal{L}_F=\emptyset$ when $n_F=0$);
\item[(3)] for any $l_1\in\mathcal{L}_F$, there is at least an $l_2\in\mathcal{L}$, $l_1\neq l_2$, such that $W^{l_1l_2}\neq0$;
hence, $\lambda_{\max}(W^{l})=1$ if $l\in\mathcal{L}_L$ and $\lambda_{\max}(W^{l})<1$ if $l\in\mathcal{L}_F$.
\end{itemize}
Naturally, when $n_F=0$ the elements in $W^{FL}$ and $W^F$ do not exist and we consider $r^{L_{n_L}}=N$.
This occurs
when all the classes are closed and recurrent and hence the state space can be partitioned into
irreducible and disjoint subspaces.
In the case of $W$ irreducible, there is only one closed and recurrent class and hence $n_L=1$ and $r^{1}=N$.


In the framework of urn systems, $W$ indicates the interacting matrix and hence
the element $w_{jh}$ represents how the color sampled from the urn $j$ is influenced by the composition of the urn $h$.
Hence, the probability of the Markov process to move from $j$ to $h$ in the state space can be interpreted as
the influence that $h$ has on $j$ in the urn system.
As a consequence, recurrent classes may be seen as subsystems of urns which are not influenced by the rest of the system;
analogously, transient classes may represent subsystems of urns which are influenced by other urns of the system.
Hence, from an interacting matrix $W$ expressed as in~\eqref{def:W},
we can decompose the urn system in:
\begin{itemize}
\item[(i)] leading systems $\{S^l,l\in\mathcal{L}_L\}$, $S^{l}:=\{r^{l^{-}}+1< j\leq r^{l}\}$,
that evolve independently with respect to the rest of the system;
\item[(ii)] if $n_F\geq0$, following systems $\{S^l,l\in\mathcal{L}_F\}$, $S^{l}:=\{r^{l^{-}}+1< j\leq r^{l}\}$, that evolve depending on the proportions of the urns in the leaders $S^{L_1},{\ldots},S^{L_{n_L}}$ and their upper followers $S^{F_1},{\ldots},S^{l^-}$.
\end{itemize}
As we will see in the following sections, the asymptotic behaviors of the leading systems
and the following systems are quite different.
For completeness of the paper, we will present the results for both the types of systems, assuming that $n_F\geq 1$.

\begin{rem}
Extensions to random and time-dependent interacting matrices $\{W_n;n\geq0\}$
are possible, as discussed in Section~\ref{section_further_extensions}.
In that case, the structure presented in~\eqref{def:W} is concerned with the limiting matrix $W:=a.s.-\lim_{n\rightarrow\infty}W_n$.
\end{rem}

\section{The interacting urn system in the stochastic approximation framework}\label{section_stochastic_approximation}

A crucial technique to characterize the asymptotic behavior of the interacting urn system consists in
revisiting its dynamics into the stochastic approximation (SA) framework.
A similar approach has been adopted in~\cite{LarPag} to establish the asymptotic behavior of a single urn.
However, since here we deal with systems of urns,
we need to extend the dynamics~\eqref{eq:dynamics_Y} to jointly study the urns that interact among each other.
To this end, we first introduce in Subsection~\ref{subection_notation_SA} a compact notation
that combines
the composition of the urns in the same subsystem $S^l$, $l\in\mathcal{L}$.
Then, in Subsection~\ref{subection_dynamics_SA} we
embed each
subsystem dynamics into the classical SA form:
given a filtered probability space $(\Omega,{\mathcal A},(\mathcal{F}_n)_{n\geq0},\mathbf{P})$, we consider the following recursive procedure
\begin{equation}\label{SAP}
\forall\, n\geq1,\quad \theta_{n}=\theta_{n-1}-\frac{1}{n}f(\theta_{n-1})+\frac{1}{n}\left(\Delta M_{n}+R_{n}\right),
\end{equation}
where $f:\mathbb{R}^d\rightarrow \mathbb{R}^d$ is a locally Lipschitz continuous function, $\theta_{n}$ an $\mathcal{F}_{n}$-measurable finite random vector and, for every $n\ge 1$,  $\Delta M_{n}$ is an $\mathcal{F}_{n-1}$-martingale increment and $R_{n}$ is an $\mathcal{F}_n$-adapted remainder term.
In our framework, the process $\theta_{n}$ satisfying~\eqref{SAP} will represent the proportions of the colors of the urns in the same subsystem.
In the next sections we
apply the ``ODE'' and the ``SDE'' methods for SA
reported in Theorem~\ref{ThmODE} and in Theorem~\ref{ThmCLT} (see Appendix),
that establish first and second-order asymptotic results for $\theta_{n}$.
Specifically, Theorem~\ref{ThmODE} states that, under suitable hypotheses on $\Delta M_{n}$ and $R_{n}$,
the set $\Theta^{\infty}$ of the limiting values of $\theta_n$ as $n\rightarrow+\infty$ is $a.s.$ a compact connected set, stable by the flow of
$ODE_f\equiv\dot{\theta}=-f(\theta)$; moreover, if $\theta^*\in\Theta^{\infty}$ is a uniformly stable equilibrium on $\Theta^{\infty}$ of $ODE_f$, then
$\theta_n\stackrel{a.s.}{\longrightarrow}\theta^*$.
In addition, under further assumptions on $\Delta M_{n}$ and $R_{n}$,
Theorem~\ref{ThmCLT} establishes the CLT for $\theta_{n}$ in which
the convergence rate and the asymptotic distribution depend on the eigen-structure of the
Jacobian matrix of $f(\theta)$ evaluated at the equilibrium point $\theta^*$.

\subsection{Notation}\label{subection_notation_SA}

The quantities related to the urn $j\in\{1,{\ldots},N\}$ at time $n$ are random variables denoted by:
\begin{itemize}
\item[(1)] $Y_{n}^{j}=(Y_{1,n}^{j},{\ldots},Y_{K,n}^{j})'\in\mathbb{R}_{+}^{K}$,
\item[(2)] $Z_{n}^{j}=(Z_{1,n}^{j},{\ldots},Z_{K,n}^{j})'\in\mathcal{S}^K$, where
$\mathcal{S}^{K}$ indicates the $K$-simplex,
\item[(3)] $\tilde{Z}_{n}^{j}=(\tilde{Z}_{1,n}^{j},{\ldots},\tilde{Z}_{K,n}^{j})'\in\mathcal{S}^K$,
\item[(4)] $X_{n}^{j}=(X_{1,n}^{j},{\ldots},X_{K,n}^{j})'\in\mathcal{S}^K\cap\{0,1\}^{K}$,
\end{itemize}
while the corresponding terms of the system $S^l$, $l\in\mathcal{L}$, given by the $s^l$ urns labeled by $\{r^{l^{-}}+1,{\ldots},r^{l}\}$,
are denoted by:
\begin{itemize}
\item[(1)] $\mathbf{Y}^{l}_n:=(Y_{n}^{r^{l^{-}}+1},{\ldots},Y_{n}^{r^{l}})'\in\mathbb{R}_{+}^{s^lK}$,
\item[(2)] $\mathbf{Z}^{l}_n:=(Z_{n}^{r^{l^{-}}+1},{\ldots},Z_{n}^{r^{l}})'\in\mathcal{S}^{s^lK}$, where
$\mathcal{S}^{s^lK}$ indicates the Cartesian product of $s^l$ $K$-simpleces where $Z_{n}^{r^{l^{-}}+1},{\ldots},Z_{n}^{r^{l}}$ are defined,
\item[(3)] $\tilde{\mathbf{Z}}^{l}_n:=(\tilde{Z}_{n}^{r^{l^{-}}+1},{\ldots},\tilde{Z}_{n}^{r^{l}})'\in\mathcal{S}^{s^lK}$,
\item[(4)] $\mathbf{X}^{l}_n:=(X_{n}^{r^{l^{-}}+1},{\ldots},X_{n}^{r^{l}})'\in\mathcal{S}^{s^lK}\cap\{0,1\}^{s^lK}$,
\item[(5)] $\mathbf{T}^{l}_n:=(T_{n}^{r^{l^{-}}+1}\mathbf{1}_K,{\ldots},T_{n}^{r^{l}}\mathbf{1}_K)'\in\mathbb{R}_{+}^{s^lK}$,
where $\mathbf{1}_K$ indicates the $K$-vector of all ones.
\end{itemize}

The replacement matrix for the system $S^l$ is defined by a non-negative
block diagonal matrix $\mathbf{D}^l_n$ of dimensions ${s^lK \times s^lK}$,
where the $s^l$ blocks are the replacement matrices of the urns $\{r^{l^{-}}+1,{\ldots},r^{l}\}$ in $S^l$, i.e. $D_{n}^{r^{l^{-}}+1},{\ldots},D_{n}^{r^{l}}$.
Analogously, the generating matrix for $S^l$ is defined by a
block diagonal matrix $\mathbf{H}^l$ of the same dimensions,
where the $s^l$ blocks are $H^{r^{l^{-}}+1},{\ldots},H_{n}^{r^{l}}$.
The interaction within the system $S^l$ is modeled by the $s^lK\times s^lK$ matrix $\mathbf{W}^l $ with values in $[0,1]$ defined as follows:
starting from $W^l$ in~\eqref{def:W}, each weight $w_{jh}$ is replaced by the corresponding diagonal matrix $w_{jh}I_K$,
where here $I_K$ indicates the $K\times K$-identity matrix.
Analogously, the interaction between a following system $S^{l_1}$, $l_1\in\mathcal{L}_F$, and another system $S^{l_2}$, $l_2\in\{L_1,{\ldots},l_1^-\}$,
is modeled by the matrix $\mathbf{W}^{l_1l_2}$,
obtained by replacing each weight $w_{jh}$ of $W^{l_1l_2}$ in~\eqref{def:W} with the corresponding diagonal matrix $w_{jh}I_K$.
Finally, we will denote by $\mathbf{I}$ the identity matrix composed by more matrices $I_K$.

\begin{Example}\label{example_1}
Consider a system of $N=2$ urns containing balls of $K=2$ colors.
Let the generating matrices $H^1$, $H^2$ and the interacting matrix $W$ be as follows:
\begin{equation}\label{def:example_H12_W}
H^1 :=
\begin{pmatrix}
3/4&1/2\\
1/4&1/2
\end{pmatrix},
\qquad
H^2 :=
\begin{pmatrix}
7/8&7/8\\
1/8&1/8
\end{pmatrix},
\qquad
W :=
\begin{pmatrix}
\alpha&1-\alpha\\
1-\beta&\beta
\end{pmatrix},
\end{equation}
where $\alpha$ and $\beta$ are given constants in $[0,1]$.

In the case of no interaction $\alpha=\beta=1$, from the classical theory on single GFUs (see~\cite{AthKar2,BaiHu,BaiHu2,Smy}), we have that
\begin{itemize}
\item[(1)] $Z^1_{n}=(Z^1_{1,n},Z^1_{2,n})'$ converges a.s. to $(2/3,1/3)'$, i.e. the right eigenvector of $H^1$ associated to $\lambda=1$; moreover
the convergence rate is $\sqrt{n}$, since the second eigenvalue of $H^1$ is $0.25$.
\item[(2)] $Z^2_{n}=(Z^2_{1,n},Z^2_{2,n})'$ converges a.s. to $(1/2,1/2)'$, i.e. the right eigenvector of $H^2$ associated to $\lambda=1$; moreover
the convergence rate is $n^{0.25}$, since the second eigenvalue of $H^2$ is $0.75$.
\end{itemize}

When both $\alpha<1$ and $\beta<1$, $W$ is irreducible.
Using the notation introduced in
Subsection~\ref{subsection_W interaction} and Subsection~\ref{subection_notation_SA},
in this case the two urns belong to the same leading subsystem $S^{L_1}=\{1,2\}$. We have
$s^{L_1}=2$,
$W=W^{L}= W^{L_1}$, and
the joint quantities read as follows:
$\mathbf{Z}^{1}_n:=(Z_{1,n}^{1},Z_{2,n}^{1},Z_{1,n}^{2},Z_{2,n}^{2})'\in\mathcal{S}^{2,2}$,
\begin{equation*}
\mathbf{H} :=
\begin{pmatrix}
\frac{3}{4} & \frac{1}{2} & 0 & 0\\ \frac{1}{4} & \frac{1}{2} & 0 & 0\\ 0 & 0 & \frac{7}{8} & \frac{1}{8}\\ 0 & 0 & \frac{1}{8} & \frac{7}{8}
\end{pmatrix}
,\qquad
\mathbf{W} :=
\begin{pmatrix}
 {\alpha} & 0 & 1 - {\alpha} & 0\\ 0 & {\alpha} & 0 & 1 - {\alpha}\\ 1 - {\beta} & 0 & {\beta} & 0\\ 0 & 1 - {\beta} & 0 & {\beta}
\end{pmatrix}.
\end{equation*}
We will discuss the asymptotic properties of this system in Example~\ref{example_2}.

When $\alpha=1$ and $\beta<1$, the first urn forms a leading system, while the second one
exhibits the behavior of a following system (see Example~\ref{example_3}).
\end{Example}

\subsection{The system dynamics in the SA form}\label{subection_dynamics_SA}

For any system $S^l$, $l\in\mathcal{L}$,
the dynamics in~\eqref{eq:dynamics_Y} can be written, using the notation of Subsection~\ref{subection_notation_SA}, as follows:
\begin{equation}\label{eq:dynamics_Y_SA}
\mathbf{Y}^l_{n}\ =\ \mathbf{Y}^l_{n-1}\ +\ \mathbf{D}^l_{n}\mathbf{X}^l_{n}.
\end{equation}
We now express~\eqref{eq:dynamics_Y_SA} in the SA form~\eqref{SAP}, where the process $\{\theta_n;n\geq1\}$ is represented by the urn proportions of the system
$S^l$, i.e. $\{\mathbf{Z}^l_{n};n\geq1\}$.
Since $\mathbf{Y}^l_{n}=diag(\mathbf{T}^l_{n})\mathbf{Z}^l_{n}$ for any $n\geq1$, from~\eqref{eq:dynamics_Y_SA} we have
$$diag(\mathbf{T}^l_{n})\mathbf{Z}^l_{n}\ =\ diag(\mathbf{T}^l_{n-1})\mathbf{Z}^l_{n-1}\ +\ \mathbf{D}^l_{n}\mathbf{X}^l_{n},$$
that is equivalent to
\begin{equation}\label{eq:dynamics_Y_SA_2}
diag(\mathbf{T}^l_{n})(\mathbf{Z}^l_{n}-\mathbf{Z}^l_{n-1})\ =\ -diag(\mathbf{T}^l_{n}-\mathbf{T}^l_{n-1})\mathbf{Z}^l_{n-1}\ +\ \mathbf{D}^l_{n}\mathbf{X}^l_{n}.
\end{equation}
Now, notice that, for any $n\geq1$,
\begin{itemize}
\item[(1)] $\mathbf{E}[diag(\mathbf{T}^l_{n}-\mathbf{T}^l_{n-1})|\mathcal{F}_{n-1}]\ =\ \mathbf{I}$ by Theorem~\ref{thm:total_number};
\item[(2)] $\mathbf{E}[\mathbf{D}^l_{n}\mathbf{X}^l_{n}|\mathcal{F}_{n-1}]\ =\ \mathbf{E}[\mathbf{D}^l_{n}|\mathcal{F}_{n-1}]\mathbf{E}[\mathbf{X}^l_{n}|\mathcal{F}_{n-1}]\ =\ \mathbf{H}^l\tilde{\mathbf{Z}}^l_{n-1}$, since $\mathbf{D}^l_{n}$ and $\mathbf{X}^l_{n}$ are independent conditionally on $\mathcal{F}_{n-1}$.
\end{itemize}
Hence, defining the martingale increment
\begin{equation}\label{eq:martingale_increment_SA}
\Delta \mathbf{M}^l_n\ :=\ \mathbf{D}^l_{n}\mathbf{X}^l_{n} - \mathbf{H}^l\tilde{\mathbf{Z}}^l_{n-1}\ -\
(diag(\mathbf{T}^l_{n}-\mathbf{T}^l_{n-1})-\mathbf{I})\mathbf{Z}^l_{n-1},
\end{equation}
we can express~\eqref{eq:dynamics_Y_SA_2} as follows:
\begin{equation}\label{eq:dynamics_Y_SA_3}
diag(\mathbf{T}^l_{n})(\mathbf{Z}^l_{n}-\mathbf{Z}^l_{n-1})\ =\ -\mathbf{Z}^l_{n-1}\ +\ \mathbf{H}^l\tilde{\mathbf{Z}}^l_{n-1}\ +\ \Delta \mathbf{M}^l_n.
\end{equation}
Now, multiplying by $diag(\mathbf{T}^l_{n})^{-1}$ and defining the remainder term
\begin{equation}\label{eq:remainder_term_SA}
\mathbf{R}^l_n\ :=\ \left(n \cdot diag(\mathbf{T}^l_{n})^{-1}-\mathbf{I}\right)
\left(-\mathbf{Z}^l_{n-1}\ +\ \mathbf{H}^l\tilde{\mathbf{Z}}^l_{n-1}\ +\ \Delta \mathbf{M}^l_n\right),
\end{equation}
we can write~\eqref{eq:dynamics_Y_SA_3} as follows:
\begin{equation}\label{eq:dynamics_Y_SA_4}
\mathbf{Z}^l_{n}-\mathbf{Z}^l_{n-1}\ =\ -\frac{1}{n}(\mathbf{Z}^l_{n-1}\ -\ \mathbf{H}^l\tilde{\mathbf{Z}}^l_{n-1})\ +\
\frac{1}{n}\left(\Delta \mathbf{M}^l_n\ +\ \mathbf{R}^l_n\right).
\end{equation}
The term $(\mathbf{Z}^l_{n-1}-\mathbf{H}^l\tilde{\mathbf{Z}}^l_{n-1})$ in~\eqref{eq:dynamics_Y_SA_4}
should represent the function $f$ in~\eqref{SAP} in the SA form.
However, although in a leader $S^l$, $l\in\mathcal{L}_L$, we have that $\tilde{\mathbf{Z}}^l_{n-1}$ only depends on $\mathbf{Z}^l_{n-1}$,
in a follower $S^l$, $l\in\mathcal{L}_F$,
the term $\tilde{\mathbf{Z}}^l_{n-1}$ is in general a function of the composition of all the urns of the system, i.e.
$\mathbf{Z}^{L_1}_{n-1},{\ldots},\mathbf{Z}^{l}_{n-1}$.
Hence, the dynamics of a leading system can be expressed as in~\eqref{eq:dynamics_Y_SA_4},
while the dynamics of a following system needs to be incorporated with other systems to be fully described.
For this reason, the asymptotic behavior of these two types of systems are studied separately:
the leading systems in Section~\ref{section_leading} and the following systems in Section~\ref{section_follower}.

\section{Leading Systems}\label{section_leading}

In this section we present the main asymptotic results concerning the leading systems $S^l$, $l\in\mathcal{L}_L$.
We recall that these systems are characterized by irreducible interacting matrices $W^l$ such that $\lambda_{max}(W^l)=1$ (see~\eqref{def:W} in Subsection~\ref{subsection_W interaction}).
For this reason, their dynamics is independent of the rest of the system and hence,
by using $\tilde{\mathbf{Z}}^l_{n-1}=\mathbf{W}^l\mathbf{Z}^l_{n-1}$ in~\eqref{eq:dynamics_Y_SA_4}, we have
\begin{equation}\label{eq:dynamics_Y_SA_leading}
\begin{aligned}
\mathbf{Z}^l_{n}-\mathbf{Z}^l_{n-1}\ &&=&\ -\frac{1}{n}h^l(\mathbf{Z}^l_{n-1})\ +\ \frac{1}{n}\left(\Delta \mathbf{M}^l_n\ +\ \mathbf{R}^l_n\right),\\
h^l(\mathbf{x})\ &&:=&\ (\mathbf{I}-\mathbf{Q}^{l})\mathbf{x},\ \ \ \mathbf{Q}^{l}\ :=\ \mathbf{H}^l\mathbf{W}^l
\end{aligned}
\end{equation}

\subsection{Extension of the urn dynamics to $\mathbb{R}^{s^l K}$}\label{subection_redefining_simplex_leading}
Since $h^l$ is defined on $\mathbb{R}^{s^l K}$,
while the process $\{\mathbf{Z}^l_{n};n\geq0\}$ takes values in the subset $\mathcal{S}^{s^l K}$,
then applying theorems based on the SA directly to~\eqref{eq:dynamics_Y_SA_leading}
may lead to improper results for the process $\mathbf{Z}^l_{n}$.
To address this issue, we appropriately modify the dynamics~\eqref{eq:dynamics_Y_SA_leading}
by replacing $h^l$ with a suitable function $f_m^l:=h^l+mg^l$,
where $m>0$ is an arbitrary constant and $g^l$ is a function defined in $\mathbb{R}^{s^l K}$ that satisfies the following properties:
\begin{itemize}
\item[(i)] the derivative $\mathcal{D}g^l$ is positive semi-definite and its kernel is $Span\{(x-y):x,y\in\mathcal{S}^{s^l K}\}$:
   hence, $g^l$ does not modify the eigen-structure of $\mathcal{D}h^l(\mathbf{x})$ on the subspace $\mathcal{S}^{s^l K}$,
   where the process $\mathbf{Z}^l_{n}$ is defined, while it changes
   the eigen-structure outside $\mathcal{S}^{s^l K}$,
   where it can be arbitrary redefined;
\item[(ii)] $g^l(\mathbf{z})=0$ for any $\mathbf{z}\in\mathcal{S}^{s^l K}$:
    hence, since $f_m^l(\mathbf{z})=h^l(\mathbf{z})$ for any $\mathbf{z}\in\mathcal{S}^{s^l K}$,
    the modified dynamics restricted to the subset $\mathcal{S}^{s^l K}$
    represents the same dynamics as in~\eqref{eq:dynamics_Y_SA_leading}.
\end{itemize}
Let us now provide an analytic expression of $g^l$.
First note that, since by definition of convex combination we always have $W^l\mathbf{1}_{s^l}=\mathbf{1}_{s^l}$,
the left eigenvectors of $W^l$ (possibly generalized) are such that $U_1'\mathbf{1}_{s^l}=1$ and $U_i'\mathbf{1}_{s^l}=0$ for all $i\neq 1$.
Denote by $Sp(A)$ the set of the eigenvalues of a matrix $A$ and note that,
since by $(\constantbalanceH)$ we always have $\mathbf{1}_K'H^j=\mathbf{1}_K'$,
then $Sp(W^l)\subset Sp(\mathbf{Q}^{l})$ and the $s^l$ left eigenvectors of $\mathbf{Q}^{l}$ associated to any $\lambda_i\in Sp(W^l)\subset Sp(\mathbf{Q}^{l})$,
$i\in\{1,{\ldots},s^l\}$, present the following structure: $\mathbf{U}_i:=(U_{i1}\mathbf{1}_K,{\ldots},U_{is^l}\mathbf{1}_K)'$.
As a consequence, for any $\mathbf{z}\in\mathcal{S}^{s^l K}$,
we have $\mathbf{U}_1'\mathbf{z}=U_1'\mathbf{1}_{s^l}=1$ and $\mathbf{U}_i'\mathbf{z}=U_i'\mathbf{1}_{s^l}=0$ for all $i\in\{2,{\ldots},s^l\}$.
Hence, denoting by $\mathbb{V}_2$ and $\mathbb{U}_2$ the matrices whose columns are $\mathbf{V}_2,{\ldots},\mathbf{V}_{s^l}$ and
$\mathbf{U}_2,{\ldots},\mathbf{U}_{s^l}$, respectively, we define the function $g^l$ as follows:
\begin{equation}\label{def:function_g}
g^l(\mathbf{x})\ :=\ \mathbf{V}_1\left(\mathbf{U}_1'\mathbf{x}-1\right)\ +\ \mathbb{V}_2\mathbb{U}_2'\mathbf{x},
\end{equation}
and the dynamics of the process $\mathbf{Z}^l_{n}$ in~\eqref{eq:dynamics_Y_SA_leading} can be replaced by the following:
\begin{equation}\label{eq:dynamics_Y_SA_leading_h}
\begin{aligned}
\mathbf{Z}^l_{n}-\mathbf{Z}^l_{n-1}\ &&=&\ -\frac{1}{n}f_m^l(\mathbf{Z}^l_{n-1})\ +\ \frac{1}{n}\left(\Delta \mathbf{M}^l_n\ +\ \mathbf{R}^l_n\right),\\
f_m^l(\mathbf{x})\ &&:=&\ (\mathbf{I}-\mathbf{Q}^{l})\mathbf{x}\ +\ m\mathbf{V}_1\left(\mathbf{U}_1'\mathbf{x}-1\right)\ +\ m\mathbb{V}_2\mathbb{U}_2'\mathbf{x}.
\end{aligned}
\end{equation}

\subsection{First-order asymptotic results}\label{subection_first_order_leading}

We now present the main convergence result concerning the limiting proportion of the urns in the leading systems.
\begin{theo}\label{thm:first_order_leading}
Assume $(\secondmomentpideltaD)$, $(\constantbalanceH)$ and $(\irreducibleH)$.
Thus, for any leading system $S^l$, $l\in\mathcal{L}_L$, we have that
\begin{equation}\label{eq:first_order_leading}
\mathbf{Z}^l_n\ \stackrel{a.s.}{\longrightarrow}\ \mathbf{Z}^l_{\infty}\ :=\ \mathbf{V}_1,
\end{equation}
where $\mathbf{V}_1$ indicates the right eigenvector associated to the simple eigenvalue $\lambda=1$ of the matrix $\mathbf{Q}^l$,
with $\sum_i V_{1i}=1$.
\end{theo}

\begin{rem}\label{rem:W_identity_Friedman_as}
Note that when the interacting matrix is the identity matrix, i.e. $W=I$, $n_L=N$ and $n_F=0$,
each urn represents a leading system and it evolves independently of the rest of the system.
In this case,~\eqref{eq:first_order_leading} expresses the usual result for a single GFU,
where the urn proportion converges to the right eigenvector associated to the maximum eigenvalue of the generating matrix, see e.g.~\cite{AthKar2,BaiHu,BaiHu2,Smy}.
\end{rem}

\begin{rem}\label{rem:simple_eigenvalue}
In Theorem~\ref{thm:first_order_leading}, condition $(\irreducibleH)$ implies that the maximum eigenvalue $\lambda=1$ of $\mathbf{Q}^l$ has multiplicity one,
which guarantees $\mathbf{V}_1$ to be the unique global attractor for the system $S^l$.
Without assumption $(\irreducibleH)$, there could be multiple attractors and hence the limiting proportions of the system would be a random variable,
as in~\cite{CriDaiPraMin,DaiPraLouisMin} where the RRU model is considered.
\end{rem}

\subsection{Second-order asymptotic results}\label{subection_second_order_leading}

We now establish the rate of convergence and the asymptotic distribution of the urn proportions in the leading systems $S^l$, $l\in\mathcal{L}_L$.
Since to obtain these results we need to apply the Central Limit Theorem of the SA (see Theorem~\ref{ThmCLT} in Appendix)
to the dynamics~\eqref{eq:dynamics_Y_SA_leading_h},
a crucial role is played by the spectrum of the $Ks^l\times Ks^l$-matrix of the first-order derivative of $f_m^l$ defined as follows:
for any $\mathbf{x}\in\mathbb{R}^{Ks^l}$
\begin{equation}\label{def:F_derivative}
\mathbf{F}^l_m\ :=\ \mathcal{D}f_m^l(\mathbf{x})\ =\ (\mathbf{I}-\mathbf{Q}^l)\ +\  m\mathbf{V}_1\mathbf{U}_1'\ +\ m\mathbb{V}_2\mathbb{U}_2'.
\end{equation}
Moreover, since the asymptotic variance depends on the second moments of the replacement matrices,
we denote by $C^j(i)$ the covariance matrix of the $i^{th}$ column of $D^j_n$, i.e.
$C^j(i):=\mathbf{C}ov[D^j_{\cdot i,n}]$, where $D^j_{\cdot i,n}:=(D^j_{1 i,n},{\ldots},D^j_{K i,n})'$;
note that $(\constantbalanceD)$ ensures the existence of $C^j(i)$.
Hence, denoting by $H^j(i):=E[H^j_{\cdot i}(H^j_{\cdot i})']$ where $H^j_{\cdot i}:=(H^j_{1 i},{\ldots},H^j_{K i})'$, we let
\begin{equation}\label{def:G}
G^j\ :=\ \sum_{i=1}^K\left(C^j(i)+H^j(i)\right)\tilde{Z}_{i,\infty}^j\ -\ Z_{\infty}^j(Z_{\infty}^j)^{'},
\end{equation}
where $\tilde{Z}_{i,\infty}^j=\sum_{h=1}^N w_{jh}Z_{i,\infty}^h$.
Then, for any leading system $S^l$, $l\in\mathcal{L}_L$,
we denote by $\mathbf{G}^l$ the block diagonal matrix made by the $s^l$ blocks $G^{r^{l^-}+1},{\ldots},G^{r^l}$.

The following theorem shows the rate of convergence and the limiting distribution of the urn proportions in the leading systems.
\begin{theo}\label{thm:second_order_leading}
Assume $(\constantbalanceD)$ and $(\irreducibleH)$.
For any leading system $S^l$, $l\in\mathcal{L}_L$, let $\lambda^{*l}$ be the eigenvalue of $Sp(\mathbf{Q}^l)\setminus Sp(W^l)$ with highest real part.
Thus, we have that $\Re e(\lambda^{*l})\equiv1-\Re e(Sp(\mathbf{F}^l_m))$ and
\begin{itemize}
\item[(a)] if $\Re e(\lambda^{*l})<1/2$, then
\begin{equation*}\label{eq:second_order_Z_less_leading}
\sqrt{n}(\mathbf{Z}^l_n-\mathbf{Z}^l_{\infty})\ \stackrel{d}{\longrightarrow}\ \mathcal{N}\left(0,\Sigma^l\right),\ \ \ \
\Sigma^l\ :=\ \lim_{m\rightarrow \infty}\int_0^{\infty}e^{u(\frac{\mathbf{I}}{2}-\mathbf{F}^l_m)}\mathbf{G}^le^{u(\frac{\mathbf{I}}{2}-\mathbf{F}^l_m)'}du.
\end{equation*}
\item[(b)] if $\Re e(\lambda^{*l})=1/2$, then
\begin{equation*}\label{eq:second_order_Z_equal_leading}
\sqrt{\frac{n}{\log(n)}}(\mathbf{Z}^l_n-\mathbf{Z}^l_{\infty})\ \stackrel{d}{\longrightarrow}\ \mathcal{N}\left(0,\Sigma^l\right).
\end{equation*}
\item[(c)] if $\Re e(\lambda^{*l})>1/2$, then there exists a finite random variable $\psi^l$ such that
\begin{equation*}\label{eq:second_order_Z_more_leading}
n^{1-\Re e(\lambda^{*l})}(\mathbf{Z}^l_n-\mathbf{Z}^l_{\infty})\ \stackrel{a.s.}{\longrightarrow}\ \psi^l.
\end{equation*}
\end{itemize}
\end{theo}

\begin{rem}\label{rem:W_identity_Friedman_clt}
When the interacting matrix $W$ is the identity matrix,
each urn represents a leading system and hence $W^l=1$ and $\mathbf{Q}^l\equiv H^l$.
In that case, $\lambda^*$ is the eigenvalue of $H^l$ with second highest real part
and hence Theorem~\ref{thm:second_order_leading} expresses the usual Central Limit Theorem for a single GFU, see e.g.~\cite{AthKar2,BaiHu,BaiHu2,Smy}.
\end{rem}

\begin{rem}\label{rem:role_W_in_Q}
The role of $\mathbf{Q}^{l}$ in Theorem~\ref{thm:second_order_leading} shows that
the convergence rate of the urns in $S^l$ does not depend only on their generating matrices $\{H^j,r^{l^-}+1\leq j\leq r^l\}$
but also on their interaction expressed in $W^l$.
For instance, consider two single GFUs whose generating matrices $H^1$ and $H^2$ are such that
the convergence rates of the urn proportions $Z_n^1$ and $Z^2_n$ without interactions are different.
Then, an interaction between these urns with an irreducible $W^l$ would make
$Z_n^1$ and $Z^2_n$ converge at the same rate,
which would depend on the choice of $W^l$.
\end{rem}

\begin{Example}[Continuation of Example~\ref{example_1}]\label{example_2}
When we introduce an interaction with an irreducible $W$,
the limit of the urn proportions changes as established in Theorem~\ref{thm:first_order_leading}.
For instance, if we consider $W$ as in~\eqref{def:example_H12_W} with $\alpha=\beta=0.8$ we have that $\mathbf{Z}_n=(Z^1_{1,n},Z^1_{2,n},Z^2_{1,n},Z^2_{2,n})'$
converges a.s. to $(0.66,0.34,0.56,0.44)'$, which is the right eigenvector of
\[
\mathbf{Q}=\mathbf{H}\mathbf{W}
=
\begin{pmatrix}
\alpha H^1 & (1-\alpha) H^1  \\
(1-\beta)H^2 & \beta H^2
\end{pmatrix}
=
\left(
\begin{smallmatrix}
\frac{3\, {\alpha}}{4} & \frac{{\alpha}}{2} & \frac{3}{4}(1 -  {\alpha}) & \frac{1}{2} ( 1- {\alpha})\\
\frac{{\alpha}}{4} & \frac{{\alpha}}{2} & \frac{1}{4}(1 -  {\alpha}) & \frac{1}{2} ( 1- {\alpha})\\
\frac{7}{8} ( 1 - {\beta}) & \frac{1}{8} ( 1 - {\beta}) & \frac{7\, {\beta}}{8} & \frac{{\beta}}{8}\\
\frac{1}{8} ( 1 - {\beta}) & \frac{7}{8} ( 1 - {\beta}) & \frac{{\beta}}{8} & \frac{7\, {\beta}}{8}
\end{smallmatrix}
\right),
\]
associated to $\lambda=1$.
Moreover, as explained in Remark~\ref{rem:role_W_in_Q}, the interaction makes the two urns
converge at the same rate, which depends on the interacting matrix, as established in Theorem~\ref{thm:second_order_leading}.
In this case $\alpha=\beta=0.8$, since $Sp(\mathbf{Q})=\{1, 0.62, 0.6, 0.18\}$ and $Sp(W)=\{1, 0.6\}$,
we have $\lambda^{*}=0.62$ and hence the convergence rate is $n^{0.38}$.
In addition, to underline the role of the interaction in the convergence rate of the system,
we note that
\begin{itemize}
\item[(i)] if $\alpha=(1-\beta)=0.8$, since $Sp(\mathbf{Q})=\{1,0.35,0,0\}$ and $Sp(W)=\{1,0\}$
we have $\lambda^{*}=0.35$ and hence the convergence rate is $\sqrt{n}$;
\item[(ii)] if $\alpha=\beta=0.5$, since $Sp(\mathbf{Q})=\{1,0.5,0,0\}$ and $Sp(W)=\{1,0\}$
we have $\lambda^{*}=0.5$ and hence the convergence rate is $\sqrt{n/\log(n)}$;
\item[(iii)] if $\alpha=(1-\beta)=0.2$, since $Sp(\mathbf{Q})=\{1,0.65,0,0\}$ and $Sp(W)=\{1,0\}$
we have $\lambda^{*}=0.65$ and hence the convergence rate is $n^{0.35}$.
\end{itemize}
\end{Example}

\section{Following Systems}\label{section_follower}

In this section we establish asymptotic properties concerning the following systems $S^l$, $l\in\mathcal{L}_F$.
As we have already underlined, the dynamics of these systems can be properly expressed in the SA form~\eqref{SAP}
only through a joint model with the urns in the systems $\{S^{L_1},{\ldots},S^l\}$.
Thus, we need a further notation to study collections of more systems.
In particular, we will replace the label $l$ with $(l)$ whenever an object is referred to the joint system
$S^{(l)}:=\{S^{L_1},{\ldots},S^l\}$ instead of the single system $S^l$.
For instance, the vector $\mathbf{Y}^{(l)}_n\in\mathbb{R}^{Kr^l}$ indicates $(\mathbf{Y}^{L_1}_n,{\ldots},\mathbf{Y}^{l}_n)'$,
and $\mathbf{D}^{(l)}_n$ indicates the block diagonal $(Kr^l\times Kr^l)$-matrix, whose blocks are made by $\mathbf{D}^{L_1}_n,{\ldots},\mathbf{D}^{l}_n$.
Then, from~\eqref{def:W} we can express the sampling probabilities in the follower $S^l$ as follows:
$$\tilde{\mathbf{Z}}^{l}_{n-1}=\sum_{i\in \{L_1,{\ldots}l^-\}} \mathbf{W}^{li}\mathbf{Z}^i_{n-1}+\mathbf{W}^{l}\mathbf{Z}^{l}_{n-1}.$$
Hence, from~\eqref{eq:dynamics_Y_SA_4} we obtain
\begin{equation}\label{eq:dynamics_Y_SA_following}
\begin{aligned}
\mathbf{Z}^{l}_{n}-\mathbf{Z}^{l}_{n-1}\ &&=&\ -\frac{1}{n}h^{l}(\mathbf{Z}^{(l^-)}_{n-1},\mathbf{Z}^{l}_{n-1})\ +
\ \frac{1}{n}\left(\Delta \mathbf{M}^{l}_n\ +\ \mathbf{R}^{l}_n\right),\\
h^{l}(\mathbf{x}_1,\mathbf{x}_2)\ &&:=&\ -\mathbf{Q}^{l(l^{-})}\mathbf{x}_1+(\mathbf{I}-\mathbf{Q}^{l})\mathbf{x}_2,\\
\mathbf{Q}^{l(l^{-})}\ &&:=&\ \left[\mathbf{H}^l\mathbf{W}^{lL_1}\ \ldots\ \mathbf{H}^l\mathbf{W}^{l\,l^-}\right],\ \
\mathbf{Q}^{l}\ :=\ \mathbf{H}^l\mathbf{W}^l
\end{aligned}
\end{equation}
Since $h^{l}$ is not only a function of $\mathbf{Z}^{l}_{n-1}$, the dynamics in~\eqref{eq:dynamics_Y_SA_following}
is not already expressed in the SA form~\eqref{SAP}.
To address this issue, we need to consider a joint model for the global system $S^{(l)}=S^{(l^-)}\cup S^l=S^{L_1}\cup{\ldots}\cup S^l$ as follows:
\begin{equation}\label{eq:dynamics_Y_SA_following_entire}\begin{aligned}
\mathbf{Z}^{(l)}_{n}-\mathbf{Z}^{(l)}_{n-1}\ &&=&\ -\frac{1}{n}h^{(l)}(\mathbf{Z}^{(l)}_{n-1})\ +\
\frac{1}{n}\left(\Delta \mathbf{M}^{(l)}_n\ +\ \mathbf{R}^{(l)}_n\right),\\
h^{(l)}(\mathbf{x})\ &&:=&\ \left(\mathbf{I}-\mathbf{Q}^{(l)}\right)\mathbf{x},
\end{aligned}\end{equation}
where $\mathbf{Q}^{(l)}$ can be recursively defined as follows:
\begin{equation}\label{def:Q_complete}
\mathbf{Q}^{(l)} :=
\begin{pmatrix}
\mathbf{Q}^{(l^{-})}&0\\
\mathbf{Q}^{l(l^{-})}&\mathbf{Q}^{l}
\end{pmatrix},
\qquad
\mathbf{Q}^{(L_{n_L})} :=
\begin{pmatrix}
\mathbf{Q}^{L_1}&{\ldots}&0\\
\ldots&\ldots&\ldots\\
0&{\ldots}&\mathbf{Q}^{L_{n_L}},
\end{pmatrix},
\end{equation}
where by convention $F_1^-=L_{n_L}$.

\subsection{Extension of the urn dynamics to $\mathbb{R}^{r^l K}$}\label{subection_redefining_simplex_follower}

We now apply to following systems similar considerations made for SA of the leading systems in Section~\ref{subection_redefining_simplex_leading}.
Note again that
$h^{(l)}$ in~\eqref{eq:dynamics_Y_SA_following_entire} is defined in $\mathbb{R}^{r^l K}$,
while the process $\{\mathbf{Z}^{(l)}_{n};n\geq0\}$ lies in the subspace $\mathcal{S}^{r^l K}$.
The application of the theorems based on the SA needs an extension of $h^{(l)}$, which takes into account the SA structure.

For this reason, we replace $h^{(l)}$ in~\eqref{eq:dynamics_Y_SA_following} with a suitable function $f_m^{(l)}:=h^{(l)}+mg^{(l)}$ such that
$m>0$ is an arbitrary constant and $g^{(l)}$ is a function defined as in~\eqref{def:function_g},
where in this case
$\{\mathbf{U}_i;1\leq i\leq r^l\}$ and $\{\mathbf{V}_i;1\leq i\leq r^l\}$,
indicate, respectively, the left and right eigenvectors of $\mathbf{Q}^{(l)}$ (possibly generalized).
Hence, the dynamics of the process $\mathbf{Z}^{(l)}_{n}$~\eqref{eq:dynamics_Y_SA_following} is replaced by the following:
\begin{equation}\label{eq:dynamics_Y_SA_following_h}\begin{aligned}
\mathbf{Z}^{(l)}_{n}-\mathbf{Z}^{(l)}_{n-1}\ &&=&\ -\frac{1}{n}f_m^{(l)}(\mathbf{Z}^{(l)}_{n-1})\ +\
\frac{1}{n}\left(\Delta \mathbf{M}^{(l)}_n\ +\ \mathbf{R}^{(l)}_n\right),\\
f_m^{(l)}(\mathbf{x})\ &&:=&\ \left(\mathbf{I}-\mathbf{Q}^{(l)}\right)\mathbf{x}\ +\ m\mathbf{V}_1\left(\mathbf{U}_1'\mathbf{x}-1\right)\ +\ m\mathbb{V}_2\mathbb{U}_2'\mathbf{x}.
\end{aligned}\end{equation}
Note that in the joint system $S^{(l)}$ the eigenvalue $\lambda=1$ of $\mathbf{Q}^{(l)}$ may not have multiplicity one;
in that case, $\mathbf{V}_1$ is univocally identified as the right eigenvector of $\mathbf{Q}^{(l)}$ associated to $\lambda=1$
such that, letting $\mathbf{U}_i:=(U_{i1}\mathbf{1}_K,{\ldots},U_{ir^l}\mathbf{1}_K)'$ and $U_i'W^{(l)}=\lambda_iU_i'$ for any $i\in\{1,{\ldots},r^l\}$,
we have $\mathbf{U}_1'\mathbf{V}_1=U_1'\mathbf{1}_{r^l}=1$ and $\mathbf{U}_i'\mathbf{V}_1=U_i'\mathbf{1}_{r^l}=0$ when $i\neq 1$.

\subsection{Removal of unnecessary components}\label{subection_in_out_follower}

The following system $S^l$ may not depend on all the components of $S^{(l^-)}$ and hence
the convergence in $S^l$ may be faster than the rate in $S^{(l^-)}$.
When this occurs, the asymptotic distribution obtained for the urn proportions in $S^{(l)}$ restricted to the urns in $S^l$ is degenerate.
To address this issue and characterize the asymptotic behavior in the following system $S^l$,
we need to reduce the dimensionality of $\mathbf{Z}^{(l)}_n$ by deleting those components which do not influence the dynamics of $\mathbf{Z}^l_n$.
Since the interaction between $S^l$ and the systems in $S^{(l^-)}$ is expressed by $\mathbf{Q^{l(l^-)}}$,
we exclude the components of $\mathbf{Z}^{(l^-)}_n$ defined on the null space of $\mathbf{Q}^{l(l^-)}$.
Formally, consider the following decomposition:
$$Sp(\mathbf{Q^{(l)}})\ =\ Sp(\mathbf{Q^{l}})\ \cup\ Sp(\mathbf{Q^{(l^-)}})\ =\ \mathcal{A}_{IN}\ \cup\ \mathcal{A}_{OUT},$$
where
\[\begin{aligned}
\mathcal{A}_{OUT}\ &&:=&\ \left\{\ \lambda\in Sp(\mathbf{Q^{(l^-)}})\ :\ \exists v \{\mathbf{Q^{(l^-)}}v=\lambda v\}\cap\{\mathbf{Q^{l(l^-)}}v=0\}\ \right\}\\
\mathcal{A}_{IN}\ &&:=&\ Sp(\mathbf{Q^{l}})\ \cup\ \left(Sp(\mathbf{Q^{(l^-)}})\setminus \mathcal{A}_{OUT}\right).
\end{aligned}\]
Then, the eigenspace of $\mathbf{Q^{(l)}}$ associated to $\lambda\in\mathcal{A}_{OUT}$ will be removed from the dynamics in~\eqref{eq:dynamics_Y_SA_following_h}.
To do this, let us denote by:
\begin{itemize}
\item[(1)] $\mathbf{U}_{IN}$ and $\mathbf{V}_{IN}$ the matrices whose columns are the left and right eigenvectors of $\mathbf{Q^{(l)}}$, respectively,
associated to eigenvalues in $\mathcal{A}_{IN}$;
\item[(2)] $\mathbf{U}_{OUT}$ and $\mathbf{V}_{OUT}$ the matrices whose columns are the left and right eigenvectors of $\mathbf{Q^{(l)}}$, respectively,
associated to eigenvalues in $\mathcal{A}_{OUT}$;
\end{itemize}
Since we do not want to modify the process $\mathbf{Z}^{(l)}_n$ on $S^l$, i.e. $\mathbf{Z}^{l}_n$,
we now construct two conjugate basis in $\mathcal{I}m(\mathbf{U}_{IN})$ and $\mathcal{I}m(\mathbf{V}_{IN})$ that are invariant on $S^l$.
Note that, since $Sp(\mathbf{Q^{l}})\subset\mathcal{A}_{IN}$,
there exists a non-singular matrix $\mathbf{P}$ such that the following decompositions hold:
\begin{equation*}
\mathbf{B}:=\mathbf{V}_{IN}\mathbf{P}=
\begin{pmatrix}
\mathbf{\hat{B}}&0\\
0&\mathbf{I}
\end{pmatrix},
\qquad
\mathbf{C}:=\mathbf{P}^{-1}\mathbf{U}_{IN}'=
\begin{pmatrix}
\mathbf{\hat{C}}&0\\
0&\mathbf{I}
\end{pmatrix}.
\end{equation*}
Since $\mathbf{\hat{C}}'\mathbf{\hat{B}}=\mathbf{I}$ and $\mathbf{\hat{B}}\mathbf{\hat{C}}'=\mathbf{V_{IN}}\mathbf{U_{IN}}'$,
$\mathbf{\hat{C}}$ and $\mathbf{\hat{B}}$ represent conjugate basis in $\mathcal{I}m(\mathbf{U}_{IN})$ and $\mathcal{I}m(\mathbf{V}_{IN})$, respectively.
Thus, for any $\mathbf{x}=(\mathbf{x}^{(l^-)},\mathbf{x}^l)'\in\mathbb{R}^{Kr^l}$, we have the following decomposition:
\begin{equation}\label{eq:decomposition_x}
\mathbf{x}\ =\ \mathbf{V_{IN}}\mathbf{U_{IN}}'\mathbf{x}\ +\ \mathbf{V_{OUT}}\mathbf{U_{OUT}}'\mathbf{x}\ =\
\mathbf{\hat{B}}\mathbf{\hat{x}}\ +\ \mathbf{V_{OUT}}\mathbf{x_{OUT}},\ \ \
\end{equation}
where
$$\mathbf{\hat{x}}:=\mathbf{\hat{C}}'\mathbf{x}=
\begin{pmatrix}
\mathbf{C}'\mathbf{x}^{(l^-)}\\
\mathbf{x}^{l}
\end{pmatrix},
\qquad \mathbf{x_{OUT}}:=\mathbf{U_{OUT}}'\mathbf{x}.$$
In particular, we consider the process $\{\mathbf{\hat{Z}}^{(l)}_{n},n\geq1\}$ defined as follows:
\begin{equation}\label{eq:decomposition_Z_hat}
\mathbf{\hat{Z}}^{(l)}_{n}:=\mathbf{\hat{C}}'\mathbf{Z}^{(l)}_{n}=
\begin{pmatrix}
\mathbf{C}'\mathbf{Z}_{n}^{(l^-)}\\
\mathbf{Z}_{n}^{l}
\end{pmatrix};
\end{equation}
now, multiplying by $\mathbf{\hat{C}}'$ to~\eqref{eq:dynamics_Y_SA_following_h}
and applying the decomposition~\eqref{eq:decomposition_x} in~\eqref{eq:dynamics_Y_SA_following_h},
since $\mathbf{\hat{C}}'\mathbb{V}_2\mathbb{U}_2'\mathbf{V_{OUT}}=0$, $\mathbf{U}_1'\mathbf{V_{OUT}}=0$ and
$\mathbf{\hat{C}}'\mathbf{V_{OUT}}=0$, we have that
\begin{equation}\label{eq:dynamics_Y_SA_following_hat}\begin{aligned}
\mathbf{\hat{Z}}^{(l)}_{n}-\mathbf{\hat{Z}}^{(l)}_{n-1}\ &&=&\ -\frac{1}{n}\hat{f}_m^{(l)}(\mathbf{\hat{Z}}^{(l)}_{n-1})\ +\
\frac{1}{n}\mathbf{\hat{C}}'\left(\Delta \mathbf{M}^{(l)}_n\ +\ \mathbf{R}^{(l)}_n\right),\\
\hat{f}_m^{(l)}(\mathbf{\hat{x}})\ &&:=&\ \left(\mathbf{I}-\mathbf{\hat{C}}'\mathbf{Q}^{(l)}\mathbf{\hat{B}}\right)\mathbf{\hat{x}}\ +\ m\mathbf{\hat{V}}_1\left(\mathbf{\hat{U}}_1'\mathbf{\hat{x}}-1\right)\ +\
m\mathbb{\hat{V}}_2\mathbb{\hat{U}}_2'\mathbf{\hat{x}},
\end{aligned}\end{equation}
where $\mathbf{\hat{U}}_1':=\mathbf{U}_1'\mathbf{\hat{B}}$, $\mathbb{\hat{U}}_2:=\mathbb{U}_2'\mathbf{\hat{B}}$, $\mathbf{\hat{V}}_1:=\mathbf{\hat{C}}'\mathbf{V}_1$ and $\mathbb{\hat{V}}_2:=\mathbf{\hat{C}}'\mathbb{V}_2$
represent the left and right eigenvectors of $\mathbf{\hat{C}}'\mathbf{Q}^{(l)}\mathbf{\hat{B}}$ associated to
$\lambda\in Sp(W^{(l)})\setminus \mathcal{A}_{OUT}$.
Since $\hat{f}_m^{(l)}$ is a function of $\mathbf{\hat{Z}}^{(l)}_{n}$,
the dynamics in~\eqref{eq:dynamics_Y_SA_following_hat} is now expressed in the SA form~\eqref{SAP}.

\begin{rem}
The interacting matrix $W$ lonely is not enough to individuate
the components of the system that actually influence a following system,
but it is necessary to study the eigen-structure of $\mathbf{Q}^{(l)}$, that joins the information of $W$ and of the generating matrices
$\{H^j,1\leq j\leq r^l\}$ of the urns in $S^{(l)}$.
This may be surprising since $W$ is the only element that defines the interaction among the urns in the system.
Nevertheless, when $H^j$ is singular, different values of $\tilde{Z}^j_{n}$ may give the same average replacements, $H^j\tilde{Z}^j_{n}$,
which is equivalent as having singularities in $W$, where different values of $\{Z^i_{n};1\leq i\leq r^l\}$
may give the same $\tilde{Z}^j_{k,n}$, and hence same $H^j\tilde{Z}^j_{n}$.
For instance, if all the columns of $H^j$ were equal to a given vector $v^j$,
the urn $j$ would be updated on average by $v^j$ regardless the value of $\tilde{Z}^j_{n-1}$
and hence the urns in $S^{(l^-)}$ would not play any role in the dynamics of the urn $j$
for any choice of $W$.
\end{rem}

\subsection{First-order asymptotic results}\label{subection_first_order_follower}

We now present the convergence result concerning the limiting proportion of the urns in the following systems.
The asymptotic behavior of $\mathbf{Z}_{n}^{(l)}$ is obtained recursively from
$\mathbf{Z}_{\infty}^{(l^-)}:=a.s.-\lim_{n\rightarrow\infty}\mathbf{Z}_{n}^{(l^-)}$.
\begin{theo}\label{thm:first_order_follower}
Assume $(\secondmomentpideltaD)$, $(\constantbalanceH)$ and $(\irreducibleH)$.
Thus, for any $l\in\mathcal{L}_F$, we have that
\begin{equation*}\label{eq:first_order_follower_joint}
\mathbf{\hat{Z}}^{(l)}_n\ \stackrel{a.s.}{\longrightarrow}\ \mathbf{\hat{Z}}^{(l)}_{\infty}\ :=\ \mathbf{\hat{V}}_1;
\end{equation*}
hence, from~\eqref{eq:decomposition_Z_hat}, in the following system $S^l$ we have that
\begin{equation*}\label{eq:first_order_follower}
\mathbf{Z}^l_n\ \stackrel{a.s.}{\longrightarrow}\ \mathbf{Z}^l_{\infty}\ :=\ \left(\mathbf{I}-\mathbf{Q}^{l}\right)^{-1}\mathbf{Q}^{l(l^-)}\mathbf{Z}^{(l^-)}_{\infty}.
\end{equation*}
\end{theo}

\subsection{Second-order asymptotic results}\label{subection_second_order_follower}

We now present the results concerning the rate of convergence and the asymptotic distribution of the urn proportions in the following systems.
To this end, let us introduce the $Ks^l\times Ks^l$-matrix of the first-order derivative of $\hat{f}_m^l$:
\begin{equation}\label{def:F_derivative_hat}
\begin{aligned}
\hat{\mathbf{F}}_m^{(l)}\ &&:=&\ \mathbf{\hat{C}}'\mathbf{F}^{(l)}_m\mathbf{\hat{B}}\\
&&=&\ (\mathbf{I}-\mathbf{\hat{C}}'\mathbf{Q}^{(l)}\mathbf{\hat{B}})\ +\  m\mathbf{\hat{V}}_1\mathbf{\hat{U}}_1'\ +\ m\mathbb{\hat{V}}_2\mathbb{\hat{U}}_2'.
\end{aligned}
\end{equation}
Moreover, the asymptotic variance will be based on the quantity $\mathbf{\hat{G}}^{(l)}:=\mathbf{\hat{C}}'\mathbf{G}^{(l)}\mathbf{\hat{B}}$,
where $\mathbf{G}^{(l)}$ is the block diagonal matrix made by $G^{1},{\ldots},G^{r^l}$ (see~\eqref{def:G}).

The following theorem shows the rate of convergence and the limiting distribution of the urn proportions in the following systems.
\begin{theo}\label{thm:second_order_follower}
Assume $(\constantbalanceD)$ and $(\irreducibleH)$.
For any following system $S^l$, $l\in\mathcal{L}_F$, let $\lambda^{*l}$ be the eigenvalue of
$Sp(\mathbf{Q}^{(l)})\setminus (Sp(W^{(l)})\cup \mathcal{A}_{OUT})$ with highest real part.
Thus, we have that $\Re e(\lambda^{*l})\equiv1-\Re e(Sp(\hat{\mathbf{F}}_m^{(l)}))$ and
\begin{itemize}
\item[(a)] if $\Re e(\lambda^{*l})<1/2$, then
\begin{equation*}\label{eq:second_order_Z_less_follower}
\sqrt{n}(\mathbf{\hat{Z}}^{(l)}_n-\mathbf{\hat{Z}}^{(l)}_{\infty})\ \stackrel{d}{\longrightarrow}\ \mathcal{N}\left(0,\hat{\Sigma}^{(l)}\right),\ \ \ \
\hat{\Sigma}^{(l)}\ :=\ \lim_{m\rightarrow \infty}\int_0^{\infty}e^{u(\frac{\mathbf{I}}{2}-\hat{\mathbf{F}}_m^{(l)})}\mathbf{\hat{G}}^{(l)}e^{u(\frac{\mathbf{I}}{2}-\hat{\mathbf{F}}_m^{(l)})'}du.
\end{equation*}
\item[(b)] if $\Re e(\lambda^{*l})=1/2$, then
\begin{equation*}\label{eq:second_order_Z_equal_follower}
\sqrt{\frac{n}{\log(n)}}(\mathbf{\hat{Z}}^{(l)}_n-\mathbf{\hat{Z}}^{(l)}_{\infty})\ \stackrel{d}{\longrightarrow}\ \mathcal{N}\left(0,\Sigma^{(l)}\right).
\end{equation*}
\item[(c)] if $\Re e(\lambda^{*l})>1/2$, then there exists a finite random variable $\psi^{(l)}$ such that
\begin{equation*}\label{eq:second_order_Z_more_follower}
n^{1-\Re e(\lambda^{*l})}(\mathbf{\hat{Z}}^{(l)}_n-\mathbf{\hat{Z}}^{(l)}_{\infty})\ \stackrel{a.s.}{\longrightarrow}\ \psi^{(l)}.
\end{equation*}
\end{itemize}
\end{theo}

\begin{rem}\label{rem:CLT_following_joint}
Note that, since from~\eqref{eq:decomposition_Z_hat} $\mathbf{\hat{Z}}^{(l)}_n=(\mathbf{C}'\mathbf{Z}^{(l^-)}_n,\mathbf{Z}^{l}_n)'$,
Theorem~\ref{thm:second_order_follower} explicitly states the limiting distribution and the asymptotic covariance structure
of the urn proportions in any following system $\mathbf{Z}^{l}_n$, $l\in\mathcal{L}_F$.
In addition, Theorem~\ref{thm:second_order_follower} also determines the correlations between $\mathbf{Z}^{l}_n$ and the components
of the urn proportions in the other systems $S^l$, $l\in\{L_1,{\ldots},l^-\}$, that actually influence the dynamics of $\mathbf{Z}^{l}_n$.
\end{rem}

\begin{rem}\label{rem:drop_irreducibility}
We highlight that condition $(\irreducibleH)$, i.e. irreducibility of the generating matrices $H^j$,
may be relaxed in Theorem~\ref{thm:first_order_follower} and Theorem~\ref{thm:second_order_follower},
by requiring $(\irreducibleH)$ only for the urns in the leading systems.
In fact, we can note from the proof that $(\irreducibleH)$
is not needed for the urns that belong to the following systems.
\end{rem}

\begin{Example}[Continuation of Example~\ref{example_1} and Example~~\ref{example_2}]\label{example_3}
Set $W$ as in~\eqref{def:example_H12_W} with $\alpha=1$ and $\beta<1$, and hence
\[
\mathbf{Q}=\mathbf{H}\mathbf{W}
=
\begin{pmatrix}
H^1 & 0 \\
(1-\beta)H^2 & \beta H^2
\end{pmatrix}.
\]
Urn 1 forms a leading system and urn 2 is a following system.
As a consequence, the asymptotic behavior of urn 1 does not depend on urn 2.
We have that
$\mathbf{Z}^{1}_{n}=(Z^1_{1,n},Z^1_{2,n})'$ converges a.s. to $\mathbf{Z}^{1}_{\infty}=(2/3,1/3)'$, and the convergence rate is $\sqrt{n}$,
see Example~\ref{example_1}.

Concerning urn 2, its limiting proportion depends also on urn 1 as established in Theorem~\ref{thm:first_order_follower},
where in this case:
\begin{equation}\label{def:example_Q_following}
\mathbf{Q}^1=H^1,\qquad\mathbf{Q}^{12}=(1-\beta)H^2,\qquad\mathbf{Q}^2=\beta H^2.
\end{equation}
For instance, if $\beta=0.5$  we have that $\mathbf{Z}^{2}_{n}=(Z^2_{1,n},Z^2_{2,n})'$
converges a.s.\ to $(\mathbf{I}-\mathbf{Q}^{2})^{-1}\mathbf{Q}^{12}\mathbf{Z}^{1}_{\infty} = (0.6,0.4)' $.
Moreover, the convergence rate of urn 2 is determined by the interaction as established in Theorem~\ref{thm:second_order_follower}.
With $\beta = 0.5$, since $Sp(\mathbf{Q})=\{1, 0.5, 0.375, 0.25\}$, $Sp(W)=\{1,0.5\}$ and $\mathcal{A}_{OUT}=\emptyset$,
we have $\lambda^{*}=0.375$ and hence the convergence rate is $\sqrt{n}$.
In addition, to underline the role of the interaction in the convergence rate of the following system,
we note that
\begin{itemize}
\item[(i)] if $\beta=0.2$, since $Sp(\mathbf{Q})=\{1,0.25,0.2,0.15\}$, $Sp(W)=\{1,0.2\}$ and $\mathcal{A}_{OUT}=\emptyset$,
we have $\lambda^{*}=0.25$ and hence the convergence rate is $\sqrt{n}$;
\item[(ii)] if $\beta=0.8$, since $Sp(\mathbf{Q})=\{1,0.8,0.6,0.25\}$ and $Sp(W)=\{1,0.8\}$ and $\mathcal{A}_{OUT}=\emptyset$,
we have $\lambda^{*}=0.6$ and hence the convergence rate is $n^{0.4}$.
\end{itemize}
If we compare these results with
the convergence rate of urn 2 without interaction ($n^{0.25}$, see Example~\ref{example_1}),
we can observe that, in this example, the interaction makes the following system converge faster.
\end{Example}

\section{Further extensions}\label{section_further_extensions}

In this section, we discuss some possible extensions of the interacting urn model presented in this paper.

\subsection{Random and time-dependent interacting matrix}

Although we consider a constant interacting matrix $W$,
the results of this paper may be extended to a system characterized by
a random sequence of interacting matrices $\{W_n;n\geq0\}$,
i.e. $W_n=[w_{jh,n}]\in\mathcal{F}_n$ and $\tilde{Z}^j_{i,n}=\sum_{h=1}^Nw_{jh,n}Z^h_{i,n}$ for any $i\in\{1,{\ldots},K\}$.
In that case, it is essential to assume the existence of a deterministic matrix $W$ such that $W_n\stackrel{a.s.}{\longrightarrow}W$,
which individuates the leading and the following systems, as in Subsection~\ref{subsection_W interaction}.

The dynamics with random and time-dependent interacting matrices could be also expressed in the SA form~\eqref{SAP},
by including the difference $(W_n-W)$ in the remainder term~\eqref{eq:remainder_term_SA}.
Naturally, the asymptotic behavior of the urn proportions would depend on the limiting interacting matrix $W$
and on the rate of convergence of the sequence $\{W_n;n\geq0\}$.
Specifically, the convergence of the urn proportions
could be obtained with the only assumption $W_n\stackrel{a.s.}{\longrightarrow}W$,
while extensions for the second-order results presented in this paper
would require $n\mathbf{E}[\|W_n-W\|^2]\rightarrow0 $ (cfr. \cite[Assumption (A5)]{LarPag}).


\subsection{Non-homogeneous generating matrices}

The independence and identically distribution of the replacement matrices is an assumption that could be relaxed
by assuming that the sequence of generating matrices $\{H_n^j;n\geq0\}$, $H^j_{n-1}:=\mathbf{E}[D^j_{n}|\mathcal{F}_{n-1}]$,
converges to some deterministic matrix $H^j$.
Thus, the urn dynamics could be expressed in the SA form~\eqref{SAP},
by including the difference $(H^j_n-H^j)$ in the remainder term~\eqref{eq:remainder_term_SA},
and the asymptotic behavior would depend on $H^j$ and on the rate of convergence of $H_n^j$.
Specifically, the second-order results would require an additional assumption as $n\mathbf{E}[\|H^j_n-H^j\|^2]\rightarrow0$
(cfr. \cite[Assumption (A5)]{LarPag}).


\section{Proofs}\label{section_proofs}

This section contains the proofs of all the results presented in the paper.
Initially, in Subsection~\ref{subsection_proofs_total_number} we prove Theorem~\ref{thm:total_number}
concerning the behavior of the total number of balls in the urns of the system.
Then, in Subsection~\ref{subsection_proofs_leading} we present the proofs of the results on the leading systems
described in Section~\ref{section_leading}.
Finally, Subsection~\ref{subsection_proofs_follower} contains the proofs of the results of Section~\ref{section_follower}
concerning the following systems.

The proofs of Subsection~\ref{subsection_proofs_leading} and~\ref{subsection_proofs_follower}
on the asymptotic behavior of the subsystems of urns are based on basic results of stochastic approximation,
which have been reported in Appendix as Theorem~\ref{ThmODE} and Theorem~\ref{ThmCLT}.

\subsection{Proof of Theorem~\ref{thm:total_number}}\label{subsection_proofs_total_number}

The proof of Theorem~\ref{thm:total_number} requires the following auxiliary result on the martingale convergence:

\begin{lem}\label{lem:preliminary}
Let $\{S_n;n\geq1\}$, $S_n:=\sum_{i=1}^n\Delta S_i$, be a zero-mean martingale with respect to a filtration $\{\mathcal{F}_n;n\geq1\}$ and
let $\{b_n;n\geq1\}$ be a non-decreasing sequence of positive numbers such that
\begin{equation}\label{eq:kronacker_assumption}
\sum_{i=1}^{\infty}b_i^{-2}\mathbf{E}[(\Delta S_i)^2|\mathcal{F}_{i-1}]\ <\ \infty,\ \ \ \ \ a.s.
\end{equation}
Then, $b_n^{-1}S_n\stackrel{a.s.}{\longrightarrow}0$.
\end{lem}

\proof
Let us define the zero-mean martingale $\tilde{S}_n:=\sum_{i=1}^n\Delta \tilde{S}_i$, with $\Delta \tilde{S}_i:=b_i^{-1}\Delta S_i$.
Equation~\eqref{eq:kronacker_assumption} states that $\sum_{i=1}^n\mathbf{E}[(\Delta \tilde{S}_i)^2|\mathcal{F}_{i-1}]<\infty$
and hence $\tilde{S}_n$ converges a.s. since its bracket $\langle\tilde{S}\rangle_{\infty}<\infty$ a.s. (see \cite[Theorem 12.13]{Will}).
Thus, the result follows by using Kronecker's Lemma (see \cite[Lemma IV.3.2]{Shi}).
\endproof

\proof[Proof of Theorem~\ref{thm:total_number}]
By using Lemma~\ref{lem:preliminary} with $b_n:=n^{1-\alpha}$ and $S_n:=T^j_{n}-n$,
the proof follows by showing that $T^j_{n}-n$ is a martingale whose increments have bounded second moments.
Now, since
$$T^j_{n}-T^j_{n-1}\ =\ \sum_{k=1}^K(Y^j_{k,n}-Y^j_{k,n-1})\ =\ \sum_{k=1}^K\sum_{i=1}^K(D^j_{ki,n} X^j_{i,n}),$$
the result follows by establishing that
\begin{itemize}
\item[(a)]\ $\sup_{n\geq1}\mathbf{E}\left[\left(\sum_{k=1}^K\sum_{i=1}^K D^j_{ki,n} X^j_{i,n}\right)^2\ \big|\mathcal{F}_{n-1}\right]<\infty$;
\item[(b)]\ $\sum_{k=1}^K\sum_{i=1}^K\mathbf{E}\left[D^j_{ki,n} X^j_{i,n}|\mathcal{F}_{n-1}\right]=1$.
\end{itemize}
For part (a), by using $|X^j_{i,n}|\leq1$ and (\secondmomentpideltaD), we have that
$$\sup_{n\geq1}\mathbf{E}\left[\left(\sum_{k=1}^K\sum_{i=1}^K (D^j_{ki,n} X^j_{i,n})\right)^2\ \big|\mathcal{F}_{n-1}\right]\ \leq\
K^2\sup_{n\geq1}\max_{j\in\{1,{\ldots},N\}}\max_{i,k\in\{1,{\ldots},K\}}\mathbf{E}\left[(D^j_{ki,n})^2\right]\ <\ \infty,$$
where the last passage follows by noticing that by Jensen's inequality and (\secondmomentpideltaD)
\begin{equation}\label{eq:relation_2_2+delta}
\mathbf{E}\left[(D^j_{ki,n})^2\right]^{\frac{1}{2}}\ \leq\ \mathbf{E}\left[(D^j_{ki,n})^{2+\delta}\right]^{\frac{1}{2+\delta}}\ <\ C_{\delta}^{\frac{1}{2+\delta}}.
\end{equation}
For part (b), since $\sum_{k=1}^KH^j_{ki}=1$ by (\constantbalanceH) and
since $D^j_{ki,n}$ and $X^j_{i,n}$ are independent conditionally on $\mathcal{F}_{n-1}$,
we obtain
$$\sum_{k=1}^K\sum_{i=1}^K\mathbf{E}\left[ D^j_{ki,n} X^j_{i,n}\ |\mathcal{F}_{n-1}\right]
=\ \sum_{k=1}^K\sum_{i=1}^K H^j_{ki} \tilde{Z}^j_{i,n-1}\
=\ \sum_{i=1}^K \tilde{Z}^j_{i,n-1} \sum_{k=1}^KH^j_{ki}\
=\ \sum_{i=1}^K \tilde{Z}^j_{i,n-1}.$$
Finally, by the definition of $\tilde{Z}^j_{i,n-1}$ in~\eqref{def:Z_tilde}, we have
$$\sum_{i=1}^K \tilde{Z}^j_{i,n-1}\ =\ \sum_{i=1}^K\sum_{h=1}^Nw_{jh}Z^{h}_{i,n-1}\ =\
\sum_{h=1}^Nw_{jh}\sum_{i=1}^KZ^{h}_{i,n-1}\ =\ \sum_{h=1}^Nw_{jh}\ =\ 1,$$
which concludes the proof of~\eqref{eq:total_number_convergence} for $\alpha<1/2$ under assumption (\constantbalanceH).

Concerning the proof of~\eqref{eq:total_number_convergence}
under assumption (\constantbalanceD), note that
$$T^j_{n}-T^j_{n-1}\ =\ \sum_{k=1}^K\sum_{i=1}^K(D^j_{ki,n} X^j_{i,n})\ =\ \sum_{i=1}^K X^j_{i,n}\ =\ 1;$$
hence, $T^j_{n}=T^j_{0}+n$ a.s. and, for any $\alpha<1$,
$$
n^{\alpha}\left(\frac{T^j_{n}}{n}-1\right)\ =\ \frac{T^j_{0}}{n^{1-\alpha}}\ \stackrel{a.s./L^2}{\longrightarrow}\ 0.\qedhere
$$
\endproof

\subsection{Proofs on the leading systems}\label{subsection_proofs_leading}

\proof[Proof of Theorem~\ref{thm:first_order_leading}]
Fix $l\in\mathcal{L}_L$ and consider the leading system $S^l=\{r^{l^-}+1\leq j\leq r^{l}\}$ with interacting matrix $W^l$.
Since the dynamic of the urn proportions $\mathbf{Z}^l_n$ in $S^l$ has been expressed in~\eqref{eq:dynamics_Y_SA_leading_h} in the SA form~\eqref{SAP},
we can establish the convergence result stated in Theorem~\ref{thm:first_order_leading} by applying Theorem~\ref{ThmODE} in Appendix.
To this end, we will show that the assumptions of Theorem~\ref{ThmODE} are satisfied by the process $\{\mathbf{Z}^l_n;n\geq1\}$ of the system $S^l$:
\begin{itemize}
\item[(1)] the function $f_m^l$ defined in~\eqref{eq:dynamics_Y_SA_leading_h} is a linear transformation and hence locally Lipschitz.
\item[(2)] from~\eqref{eq:martingale_increment_SA}, we have that
$\sup_{n\geq 1}\mathbf{E}\left[\left\|\Delta \mathbf{M}^l_n\right\|^2|\mathcal{F}_{n-1}\right]<\infty$ is satisfied by establishing
\begin{itemize}
\item[(2a)] $\sup_{n\geq 1}\mathbf{E}\left[\left\|\mathbf{D}^l_{n}\mathbf{X}^l_{n}\right\|^2|\mathcal{F}_{n-1}\right]<\infty$;
\item[(2b)] $\sup_{n\geq 1}\mathbf{E}\left[\left\|diag(\mathbf{T}^l_{n}-\mathbf{T}^l_{n-1})\mathbf{Z}^l_{n-1}\right\|^2
|\mathcal{F}_{n-1}\right]<\infty$.
\end{itemize}
Concerning (2a), since $X^j_{i,n}\in\{0,1\}$ a.s., we have that
$$\left\|\mathbf{D}^l_{n}\mathbf{X}^l_{n}\right\|^2\ \leq\ \sum_{j\in S^l}\sum_{k=1}^K\sum_{i=1}^K\left(D^j_{ki,n}\right)^2,\ \ a.s.$$
Thus, (2a) follows by assumption (\secondmomentpideltaD), since
$$\sup_{n\geq 1}\mathbf{E}\left[\left\|\mathbf{D}^l_{n}\mathbf{X}^l_{n}\right\|^2|\mathcal{F}_{n-1}\right]\ \leq\
\sum_{j\in S^l}\sum_{k=1}^K\sum_{i=1}^K\sup_{n\geq 1}\mathbf{E}\left[\left(D^j_{ki,n}\right)^2\right]\ \leq\ s^lK^2C_{\delta}^{\frac{2}{2+\delta}},$$
where $\mathbf{E}[(D^j_{ki,n})^2]\leq C_{\delta}^{\frac{2}{2+\delta}}$ follows by~\eqref{eq:relation_2_2+delta}.
Concerning (2b), since $\sum_{i=1}^K(Z^j_{i,n})^2\leq1$, we have
\begin{equation}\label{eq:step1_2b}
\left\|diag(\mathbf{T}^l_{n}-\mathbf{T}^l_{n-1})\mathbf{Z}^l_{n-1}\right\|^2\ \leq\ \sum_{j\in S^l}(T^j_n-T^j_{n-1})^2,\ \ a.s.
\end{equation}
where we recall that
\begin{equation}\label{eq:step2_2b}
T^j_{n}-T^j_{n-1}\ =\ \sum_{k=1}^K(Y^j_{k,n}-Y^j_{k,n-1})\ =\ \sum_{k=1}^K\sum_{i=1}^K(D^j_{ki,n} X^j_{i,n}).
\end{equation}
Hence, combining~\eqref{eq:step1_2b} and~\eqref{eq:step2_2b}, since $X^j_{i,n}\in\{0,1\}$ and $\sum_{i=1}^KX^j_{i,n}=1$ a.s., we obtain that
\begin{equation*}\begin{aligned}
\left\|diag(\mathbf{T}^l_{n}-\mathbf{T}^l_{n-1})\mathbf{Z}^l_{n-1}\right\|^2\
&&\leq&\ \sum_{j\in S^l}\left(\sum_{k=1}^K\sum_{i=1}^K (D^j_{ki,n}X^j_{i,n})\right)^2\\
&&\leq&\ \sum_{j\in S^l}\sum_{i=1}^K\left(\sum_{k=1}^K D^j_{ki,n}\right)^2,\ \ a.s.
\end{aligned}\end{equation*}
Finally, using the relation $(\sum_{k=1}^K a_k^2)\leq K^2(\sum_{k=1}^Ka_k^2)$, (2b) follows by assumption (\secondmomentpideltaD), since
\begin{multline*}
\sup_{n\geq 1}\mathbf{E}\left[\left\|diag(\mathbf{T}^l_{n}-\mathbf{T}^l_{n-1})\mathbf{Z}^l_{n-1}\right\|^2|\mathcal{F}_{n-1}\right]\\
\leq
\sup_{n\geq 1}\sum_{j\in S^l}\sum_{i=1}^KK^2\sum_{k=1}^K\mathbf{E}\left[\left(D^j_{ki,n}\right)^2\right]\ \leq\ s^lK^4 C_{\delta}^{\frac{2}{2+\delta}},
\end{multline*}
where $\mathbf{E}[(D^j_{ki,n})^2]\leq C_{\delta}^{\frac{2}{2+\delta}}$ follows by~\eqref{eq:relation_2_2+delta}.
\item[(3)] from~\eqref{eq:remainder_term_SA}, we show $\|\mathbf{R}^l_n\|\stackrel{a.s.}{\longrightarrow}0$ by establishing that,
for any $(2+\delta)^{-1}<\alpha<2^{-1}$,
\begin{itemize}
\item[(3a)] $n^{\alpha}\left\|n\cdot diag(\mathbf{T}^l_{n})^{-1}-\mathbf{I}\right\|\stackrel{a.s.}{\longrightarrow}0$,
\item[(3b)] $n^{-\alpha}\left\|\mathbf{Z}^l_{n-1}-\mathbf{H}^l\tilde{\mathbf{Z}}^l_{n-1}\right\|\stackrel{a.s.}{\longrightarrow}0$,
\item[(3c)] $n^{-\alpha}\left\|\Delta \mathbf{M}^l_n\right\|\stackrel{a.s.}{\longrightarrow}0$,
\end{itemize}
where we recall that $\delta>0$ is defined in Assumption (\secondmomentpideltaD) (see Subsection~\ref{subection_main_assumptions}).
Since (3a) follows straightforwardly by Theorem~\ref{thm:total_number}, consider (3b).
For any $\epsilon>0$, using Markov's inequality we obtain
\begin{equation*}
\mathbf{P}\left(\left\|\mathbf{Z}^l_{n-1}-\mathbf{H}^l\tilde{\mathbf{Z}}^l_{n-1}\right\|>\epsilon n^{\alpha}\right)\ \leq\
(\epsilon n^{\alpha})^{-(2+\delta)}\mathbf{E}\left[\left\|\mathbf{Z}^l_{n-1}-\mathbf{H}^l\tilde{\mathbf{Z}}^l_{n-1}\right\|^{(2+\delta)}\right].
\end{equation*}
Hence, (3b) follows by Borel-Cantelli Lemma since $\alpha\cdot(2+\delta)>1$ and
$$\sup_{n\geq0}\mathbf{E}\left[\left\|\mathbf{Z}^l_{n-1}-\mathbf{H}^l\tilde{\mathbf{Z}}^l_{n-1}\right\|^{(2+\delta)}\right]\
\leq\ \sum_{j\in S^l}2^{(2+\delta)}\ <\ \infty.$$
Concerning (3c), we can apply again Markov's inequality and the same arguments of part (3b) since by assumption (\secondmomentpideltaD) we have that
$$\sup_{n\geq0}\mathbf{E}\left[\left\|\mathbf{D}^l_{n}\mathbf{X}^l_{n} - \mathbf{H}^l\tilde{\mathbf{Z}}^l_{n-1}\right\|^{(2+\delta)}\right]\ \leq\
\sup_{n\geq0}\sum_{j\in S^l}\sum_{k=1}^K\sum_{i=1}^K\mathbf{E}\left[(D^l_{ki,n})^{(2+\delta)}\right] \ <\ \infty.$$
\end{itemize}

Thus, by applying Theorem~\ref{ThmODE} to the dynamics in~\eqref{eq:dynamics_Y_SA_leading_h},
we have that the limiting values of $\mathbf{Z}^l_n$ are included in the set
$$\left\{\ \mathbf{x}\in\mathbb{R}^{Ks_l}\ :\ f_m^l(\mathbf{x})=0\ \right\}.$$
Now, denote by $\mathbb{V}_3$ and $\mathbb{U}_3$ the matrices whose columns are, respectively, the right and left eigenvectors of $\mathbf{Q}^l$
(possibly generalized) associated to the eigenvalues $\lambda\in Sp(\mathbf{Q}^l)\setminus Sp(W^l)$.
Hence, we have the following decomposition
\begin{equation}\label{eq:decomposition_Q}
\mathbf{Q}^l\ =\ \mathbf{V}_1\mathbf{U}_1'\ +\ \mathbb{V}_2\mathbf{J}_2\mathbb{U}_2'\ +\ \mathbb{V}_3\mathbf{J}_3\mathbb{U}_3',
\end{equation}
where $\mathbf{J}_2$ and $\mathbf{J}_3$ represent the corresponding jordan blocks.
Since the eigenvectors of $\mathbf{Q}^l$ represent a basis of $\mathbb{R}^{Ks_l}$,
for any $\mathbf{x}\in\mathbb{R}^{Ks_l}$ there exists $a\in \mathbb{R}$, $b\in \mathbb{R}^{s^l-1}$ and $c\in \mathbb{R}^{s^l(K-1)}$ such that
\begin{equation}\label{eq:decomposition_x_initial}
\mathbf{x}\ =\ \mathbf{V}_1a\ +\ \mathbb{V}_2b\ +\ \mathbb{V}_3c.
\end{equation}
Hence, by using~\eqref{eq:decomposition_Q} and~\eqref{eq:decomposition_x_initial}, we obtain
\[\begin{aligned}
h^l(\mathbf{x})\ &&=&\ \mathbb{V}_2(\mathbf{I}-\mathbf{J}_2)b\ +\ \mathbb{V}_3(\mathbf{I}-\mathbf{J}_3)c,\\
g^l(\mathbf{x})\ &&=&\ \mathbf{V}_1(a-1)\ +\ \mathbb{V}_2b,
\end{aligned}\]
and then, since $f_m^l(\mathbf{x})=h^l(\mathbf{x})+mg^l(\mathbf{x})$, it gives us
\begin{equation}\label{eq:h_transformed_x}
f_m^l(\mathbf{x})\ =\ m\mathbf{V}_1(a-1)\ +\ \mathbb{V}_2((1+m)\mathbf{I}-\mathbf{J}_2)b\ +\ \mathbb{V}_3(\mathbf{I}-\mathbf{J}_3)c.
\end{equation}
From the irreducibility of $H^j$ assumed in (\irreducibleH), for all $\lambda\in Sp(\mathbf{Q}^l)\setminus Sp(W^l)$ we have $\lambda<1$
and hence $(\mathbf{I}-\mathbf{J}_3)$ is positive definite.
Therefore, since $m>0$, from~\eqref{eq:h_transformed_x} we have that $f_m^l(\mathbf{x})=0$ if and only if $a=1$ and $b=c=0$, i.e. $\mathbf{x}=\mathbf{V}_1$.

It
remains to prove that $\mathbf{V}_1$ is a global attractor in $\mathbb{R}^{Ks_l}$.
To this end, we will show that the Jacobian matrix
$\mathcal{D}f_m^l(\mathbf{x})$ is positive definite for any $\mathbf{x}\in\mathbb{R}^{Ks_l}$.
We recall that, from~\eqref{def:F_derivative} we have
\begin{equation}\label{eq:F_derivative_jordan}
\mathbf{F}^l_m\ =\ \mathcal{D}f_m^l(\mathbf{x})\ =\
m\mathbf{V}_1\mathbf{U}_1'\ +\ \mathbb{V}_2((1+m)\mathbf{I}-\mathbf{J}_2)\mathbb{U}_2'\ +\ \mathbb{V}_3(\mathbf{I}-\mathbf{J}_3)\mathbb{U}_3'.
\end{equation}
Hence, since $m>0$ and $(\mathbf{I}-\mathbf{J}_3)$ is positive definite by assumption (\irreducibleH),
we have that $\mathbf{F}^l_m$ is positive definite for any $m>0$.
This concludes the proof.
\endproof

\proof[Proof of Theorem~\ref{thm:second_order_leading}]
The proof consists in showing that the assumptions of the $CLT$ for processes in the $SA$ form (Theorem~\ref{ThmCLT} in Appendix) are satisfied by
the dynamics in~\eqref{eq:dynamics_Y_SA_following_h} of the urn proportions $\mathbf{Z}^l_{n}$ in the leading system $S^l$.\\

First, we show that condition $\{\Re e(Sp(Df(\theta^*)))>1/2\}$ in Theorem~\ref{ThmCLT} is equivalent to $\{\Re e(\lambda^{*l})<1/2\}$.
Note that the function $f$ of the SA form~\eqref{SAP} is represented in our case by $f_m^l$ defined in~\eqref{eq:dynamics_Y_SA_leading_h}.
Similarly, the term $\theta^*$ in Appendix indicates the deterministic limiting proportion $\mathbf{Z}^l_{\infty}$,
while $Dh(\theta^*)$ is represented by $\mathbf{F}^l_m$ defined in~\eqref{def:F_derivative}.

Now, consider the eigen-structure of $\mathbf{Q}^l$ and note that
$\mathbf{F}^l_m$ has been expressed in~\eqref{eq:F_derivative_jordan} as follows:
$$\mathbf{F}^l_m\ =\ m\mathbf{V}_1\mathbf{U}_1'\ +\ \mathbb{V}_2((1+m)\mathbf{I}-\mathbf{J}_2)\mathbb{U}_2'\ +\ \mathbb{V}_3(\mathbf{I}-\mathbf{J}_3)\mathbb{U}_3',$$
Hence, it is easy to see that the eigenvectors of $\mathbf{F}^l_m$ and $\mathbf{Q}^l$ are the same, since
\begin{itemize}
\item[(1)] $\mathbf{F}^l_m\mathbf{V}_1=m\mathbf{V}_1$,
\item[(2)] $\mathbf{F}^l_m\mathbb{V}_2=\mathbb{V}_2((1+m)\mathbf{I}-\mathbf{J}_2)$,
\item[(3)] $\mathbf{F}^l_m\mathbb{V}_3=\mathbb{V}_3(\mathbf{I}-\mathbf{J}_3)$.
\end{itemize}
Thus
$$Sp(\mathbf{F}^l_m)=\{m\}\cup \left\{(1+m)-\lambda,\lambda\in Sp(W^l)\setminus\{1\}\right\}\cup \left\{1-\lambda,\lambda\in Sp(\mathbf{Q}^l)\setminus Sp(W^l)\right\}.$$
By setting $m>0$ arbitrary large,
we obtain that
\[
\{\Re e(Sp(\mathcal{D}f(\theta^*)))>1/2\}\equiv\{\Re e(\lambda^{*l})<1/2\}.
\]
Condition~\eqref{HypDM1} of Theorem~\ref{ThmCLT} follows from
analogous arguments of point (2) in the proof of Theorem~\ref{thm:first_order_leading}.
In fact, since
$$\sup_{n\geq1}\mathbf{E}[\|\Delta\mathbf{M}^l_{n}\|^{2+\delta}|\mathcal{F}_{n-1}]\ \leq\
K^{2+\delta}\sum_{j=1}^N\sum_{i=1}^K\sum_{k=1}^K\sup_{n\geq1}\mathbf{E}[(D^j_{ki,n})^{2+\delta}]
\ \leq\ NK^{4+\delta}.$$
For what concerns condition~\eqref{HypDM2}, we will show in a moment that for any $l\in\mathcal{L}_L$
$$\mathbf{E}[\Delta\mathbf{M}^l_{n}(\Delta\mathbf{M}^l_{n})^{'}|\mathcal{F}_{n-1}]\stackrel{a.s.}{\longrightarrow}\mathbf{G}^l,\ \quad\
\mathbf{E}[\Delta\mathbf{M}^{l_1}_{n}(\Delta\mathbf{M}^{l_2}_{n})^{'}]=0\ \ \ \forall l_1\neq l_2.$$
To this end, we first show that, for any urn $j\in S^l$,
$\mathbf{E}[\Delta M_{n}^j(\Delta M_{n}^j)'|\mathcal{F}_{n-1}] \stackrel{a.s.}{\longrightarrow} G^j$.
Note that
$$
\mathbf{E}[\Delta M_{n}^j(\Delta M_{n}^j)^{'}|\mathcal{F}_{n-1}]\ =\ \mathbf{E}[(D_n^j X_{n}^j)(D_n^jX_{n}^j)^{'}|\mathcal{F}_{n-1}]\ -\
(H^j\tilde{Z}_{n-1}^j)(H^j\tilde{Z}_{n-1}^j)^{'},
$$
and the first term of the right-hand side can be written as
\[\begin{aligned}\mathbf{E}[(D_n^j X_{n}^j)(D_n^j X_{n}^j)^{'}|\mathcal{F}_{n-1}]\
&&=&\ \sum_{i=1}^K\mathbf{E}[D_{\cdot i,n}^j(D_{\cdot i,n}^j)^{'}|\mathcal{F}_{n-1}]\mathbf{P}(X_{i,n}^j=1|\mathcal{F}_{n-1})\\
&&=&\ \sum_{i=1}^K(C^j(i)+H^j(i))\tilde{Z}_{i,n}^j.
\end{aligned}\]
When
$n$ increases to infinity, from~\eqref{def:G} we obtain
$$\mathbf{E}[\Delta M_{n}^j(\Delta M_{n}^j)^{'}|\mathcal{F}_{n-1}]\ \stackrel{a.s.}{\longrightarrow}\
\sum_{i=1}^K(C^j(i)+H^j(i))\tilde{Z}_{i,\infty}^j\ -\
Z_{\infty}^j(Z_{\infty}^j)^{'}\ =\ G^j.$$
We recall that for any $j_1\neq j_2$, $D_n^{j_1} X_{n}^{j_1}$ and $D_n^{j_2} X_{n}^{j_2}$ are independent conditionally on $\mathcal{F}_{n-1}$.
As a consequence,
$\mathbf{E}[\Delta M_{n}^{j_1}(\Delta M_{n}^{j_2})^{'}|\mathcal{F}_{n-1}]=0$ and hence
$\mathbf{E}[\Delta\mathbf{M}^{l_1}_{n}(\Delta\mathbf{M}^{l_2}_{n})^{'}]=0$ for any $l_1\neq l_2$.\\

It remains to check that the remainder sequence $\{\mathbf{R}^l_n;n\geq1\}$ satisfies~(\ref{HypReste}) for any $\epsilon>0$, i.e.
\begin{equation}\label{eq:check_remainder_CLT}
\mathbf{E}\left[n\|\mathbf{R}^l_{n}\|^2\mathds{1}_{\{\left\|\mathbf{Z}^l_n-\mathbf{Z}^l_{\infty}\right\|\leq \epsilon\}}\right]\longrightarrow 0.
\end{equation}
Equation \eqref{eq:check_remainder_CLT} can be obtained by
combining~\eqref{eq:remainder_term_SA} and part (3b) in the proof of Theorem~\ref{thm:first_order_leading},
once we have observed that assumption (\constantbalanceD) in Theorem~\ref{thm:total_number} implies that
\begin{equation*}
\mathbf{E}\left[n\left\|n\cdot diag(\mathbf{T}^l_{n})^{-1}-\mathbf{I}\right\|^2\right]\longrightarrow0 .
\end{equation*}

Since the assumptions are all satisfied, we can apply Theorem~\ref{ThmCLT} to any leading system $S^l$, $l\in\mathcal{L}_L$,
so obtaining the CLT of Theorem~\ref{thm:second_order_leading}, with asymptotic variance
$$\Sigma_m^l\ :=\
\int_0^{\infty}e^{u(\frac{\mathbf{I}}{2}-\mathbf{F}^l_m)}\mathbf{G}^le^{u(\frac{\mathbf{I}}{2}-\mathbf{F}^l_m)'}du.$$
Finally, we need to fix $m>0$ to obtain the correct asymptotic variance $\Sigma^l$ for
$(\mathbf{Z}^l_n-\mathbf{Z}^l_{\infty})$ in $Span\{(x-y):x,y\in\mathcal{S}^{s^l K}\}$.
Since by construction the kernel of $\mathbf{U}'_i$, $i\in\{r^{l^-}+1,r^l\}$, is exactly $Span\{(x-y):x,y\in\mathcal{S}^{s^l K}\}$,
we impose that $\mathbf{U}'_i\Sigma_m^l\mathbf{U}_i=0$ so obtaining that
$\Sigma^l=\lim_{m\rightarrow\infty}\Sigma_m^l$.
This concludes the proof.
\endproof

\subsection{Proofs on the following systems}\label{subsection_proofs_follower}

\proof[Proof of Theorem~\ref{thm:first_order_follower}]
Consider the joint system $S^{(l)}=\cup_{i\in\{L_1,{\ldots}l\}}S^i$, for $l\in\mathcal{L}_F$,
composed by the leading systems $S^{L_1},{\ldots}S^{L_{n_L}}$ and the following systems $S^{F_1},{\ldots}S^{l}$, where we recall $S^{l}:=\{r^{l^-}+1\leq j\leq r^l\}$.
As explained in Section~\ref{section_follower}, we focus on the reduced process $\mathbf{\hat{Z}}_{n}^{(l)}:=\mathbf{\hat{C}}'\mathbf{Z}_{n}^{(l)}$,
whose dynamics is expressed in~\eqref{eq:dynamics_Y_SA_following_hat} as follows:
\begin{equation}\label{eq:dynamics_Y_SA_following_hat_proof_first}\begin{aligned}
\mathbf{\hat{Z}}^{(l)}_{n}-\mathbf{\hat{Z}}^{(l)}_{n-1}\ &&=&\ -\frac{1}{n}\hat{f}_m^{(l)}(\mathbf{\hat{Z}}^{(l)}_{n-1})\ +\
\frac{1}{n}\mathbf{\hat{C}}'\left(\Delta \mathbf{M}^{(l)}_n\ +\ \mathbf{R}^{(l)}_n\right),\\
\hat{f}_m^{(l)}(\mathbf{\hat{x}})\ &&:=&\ \left(\mathbf{I}-\mathbf{\hat{C}}'\mathbf{Q}^{(l)}\mathbf{\hat{B}}\right)\mathbf{\hat{x}}\ +\ m\mathbf{\hat{V}}_1\left(\mathbf{\hat{U}}_1'\mathbf{\hat{x}}-1\right)\ +\
m\mathbb{\hat{V}}_2\mathbb{\hat{U}}_2'\mathbf{\hat{x}},
\end{aligned}\end{equation}
where the function $f$ in the SA form~\eqref{SAP} is here represented by $\hat{f}_m^{(l)}$ that takes values in $Span\{\mathbf{V}_{IN}\}$.

Analogously to the proof of Theorem~\ref{thm:first_order_leading} for the leading systems, one can show that
all the assumptions of Theorem~\ref{ThmODE} are satisfied by the dynamics in~\eqref{eq:dynamics_Y_SA_following_hat_proof_first} and hence
the limiting values of $\mathbf{\hat{Z}}^{(l)}_{n}$ are represented by those $\mathbf{x}\in Span\{\mathbf{V}_{IN}\}$ such that $\hat{f}_m^{(l)}(\mathbf{x})=0$.
We use analogous decompositions of those in~\eqref{eq:decomposition_Q}
and in~\eqref{eq:decomposition_x_initial} for $\mathbf{\hat{C}}'\mathbf{Q}^{(l)}\mathbf{\hat{B}}$
and $\mathbf{x}\in Span\{\mathbf{V}_{IN}\}$ respectively, obtaining
\begin{equation}\label{eq:h_transformed_x_follower}
\hat{f}_m^l(\mathbf{x})\ =\ m\mathbf{\hat{V}}_1(a-1)\ +\ \mathbb{\hat{V}}_2((1+m)\mathbf{I}-\mathbf{\hat{J}}_2)b\ +\ \mathbb{\hat{V}}_3(\mathbf{I}-\mathbf{\hat{J}}_3)c,
\end{equation}
where $\mathbf{\hat{J}}_2:=\mathbf{\hat{C}}'\mathbf{J}_2\mathbf{\hat{B}}$ and $\mathbf{\hat{J}}_3:=\mathbf{\hat{C}}'\mathbf{J}_3\mathbf{\hat{B}}$.
By assumption (\irreducibleH), $H^j$ are irreducible.
Thus, $\lambda<1$ for all $\lambda\in \mathcal{A}_{IN}\setminus Sp(W^{(l)})$
and hence $(\mathbf{I}-\mathbf{\hat{J}}_3)$ is positive definite.
Therefore, since $m>0$, from~\eqref{eq:h_transformed_x_follower} we have that $\hat{f}_m^l(\mathbf{x})=0$ if and only if $a=1$ and $b=c=0$, i.e. $\mathbf{x}=\mathbf{\hat{V}}_1$.


By definition of $\mathbf{Q}^{(l)}$ (see \eqref{def:Q_complete}), we have that
$$
\mathbf{\hat{C}}'\mathbf{Q}^{(l)}\mathbf{\hat{B}} =
\begin{pmatrix}
\mathbf{C}'\mathbf{Q}^{(l^{-})}\mathbf{B}&0\\
\mathbf{Q}^{l(l^{-})}\mathbf{B}&\mathbf{Q}^{l}
\end{pmatrix},
$$
and hence we can express $\mathbf{\hat{V}}_1=(\mathbf{\hat{V}}^{(l^-)}_1,\mathbf{\hat{V}}^{l}_1)'$ as follows:
$$
\mathbf{\hat{V}}_1 =
\begin{pmatrix}
\mathbf{\hat{V}}^{(l^-)}_1\\
(\mathbf{I}-\mathbf{Q}^{l})^{-1}\mathbf{Q}^{l(l^-)}\mathbf{B}\mathbf{\hat{V}}^{(l^-)}_1
\end{pmatrix}
=
\begin{pmatrix}
\mathbf{C}'\mathbf{V}^{(l^-)}_1\\
(\mathbf{I}-\mathbf{Q}^{l})^{-1}\mathbf{Q}^{l(l^-)}\mathbf{B}\mathbf{C}'\mathbf{V}^{(l^-)}_1
\end{pmatrix}.
$$
Now, since $\mathbf{V}^{(l^-)}_1\in \mathcal{I}m(\mathbf{V}_{IN})$, we have
\[
\mathbf{B}\mathbf{C}'\mathbf{V}^{(l^-)}_1=\mathbf{V}_{IN}\mathbf{U}_{IN}'\mathbf{V}^{(l^-)}_1=\mathbf{V}^{(l^-)}_1.
\]

Finally, since from~\eqref{def:F_derivative_hat} $\hat{\mathbf{F}}_m^{(l)}=\mathbf{\hat{C}}'\mathbf{F}^{(l)}_m\mathbf{\hat{B}}$,
we have that $Sp(\hat{\mathbf{F}}^{(l)}_m)\subset Sp(\mathbf{F}^{(l)}_m)$
and hence $\hat{\mathbf{F}}_m^{(l)}$ is positive definite for any $m>0$.
As a consequence, $\mathbf{\hat{V}}_1$ is a global attractor in $Span\{\mathbf{V}_{IN}\}$ and
this concludes the proof.
\endproof

\begin{rem}
We highlight that, when (\irreducibleH) does not hold, the matrix $(\mathbf{I}-\mathbf{J}_3)$ in~\eqref{eq:h_transformed_x_follower}
may not be positive definite and hence the solution $\mathbf{\hat{V}}_1$ would not be unique.
However, since in the following systems $S^l$, $l\in\mathcal{L}_F$, we have $\lambda_{\max}(W^{l})<1$ and this implies $\lambda_{\max}(\mathbf{Q}^{l})<1$,
the irreducibility assumption of $H^j$ required in (\irreducibleH) is not necessary for the following systems,
but it is only essential in the leading systems in which $\lambda_{\max}(W^{l})=1$.
\end{rem}

\proof[Proof of Theorem~\ref{thm:second_order_follower}]
Consider the joint system $S^{(l)}=\cup_{i\in\{L_1,{\ldots}l\}}S^i$, for $l\in\mathcal{L}_F$,
composed by the leading systems $S^{L_1},{\ldots}S^{L_{n_L}}$ and the following systems $S^{F_1},{\ldots}S^{l}$, where we recall $S^{l}:=\{r^{l^-}+1\leq j\leq r^l\}$.
As explained in Section~\ref{section_follower}, we focus on the reduced process $\mathbf{\hat{Z}}_{n}^{(l)}:=\mathbf{\hat{C}}'\mathbf{Z}_{n}^{(l)}$,
whose dynamics is expressed in~\eqref{eq:dynamics_Y_SA_following_hat} as follows:
\[\begin{aligned}
\mathbf{\hat{Z}}^{(l)}_{n}-\mathbf{\hat{Z}}^{(l)}_{n-1}\ &&=&\ -\frac{1}{n}\hat{f}_m^{(l)}(\mathbf{\hat{Z}}^{(l)}_{n-1})\ +\
\frac{1}{n}\mathbf{\hat{C}}'\left(\Delta \mathbf{M}^{(l)}_n\ +\ \mathbf{R}^{(l)}_n\right),\\
\hat{f}_m^{(l)}(\mathbf{\hat{x}})\ &&:=&\ \left(\mathbf{I}-\mathbf{\hat{C}}'\mathbf{Q}^{(l)}\mathbf{\hat{B}}\right)\mathbf{\hat{x}}\ +\ m\mathbf{\hat{V}}_1\left(\mathbf{\hat{U}}_1'\mathbf{\hat{x}}-1\right)\ +\
m\mathbb{\hat{V}}_2\mathbb{\hat{U}}_2'\mathbf{\hat{x}},
\end{aligned}\]
where the function $f$ in the SA form in~\eqref{SAP} is here represented by $\hat{f}_m^{(l)}$.
The proof will be realized by showing that the assumptions of the Theorem~\ref{ThmCLT} in Appendix are satisfied by
the process $\mathbf{\hat{Z}}_{n}^{(l)}$, with $\theta^*$ replaced by the deterministic limiting proportion $\mathbf{\hat{Z}}^{(l)}_{\infty}$,
and $\mathcal{D}f(\theta^*)$ represented by $\hat{\mathbf{F}}_m^{(l)}$ defined in~\eqref{def:F_derivative_hat}.

To do this,  we first show that condition $\{\Re e(Sp(Df(\theta^*)))>1/2\}$ in Theorem~\ref{ThmCLT} is equivalent to $\{\Re e(\lambda^{*l})<1/2\}$.
To this end, analogously to the proof of Theorem~\ref{thm:second_order_leading} for the leading systems,
note that
\begin{itemize}
\item[(1)] $\hat{\mathbf{F}}_m^{(l)}\mathbf{\hat{V}}_1=m\mathbf{\hat{V}}_1$,
\item[(2)] $\hat{\mathbf{F}}_m^{(l)}\mathbb{\hat{V}}_2=\mathbb{\hat{V}}_2((1+m)\mathbf{I}-\mathbf{J}_2)$,
\item[(3)] $\hat{\mathbf{F}}_m^{(l)}\mathbb{\hat{V}}_3=\mathbb{\hat{V}}_3(\mathbf{I}-\mathbf{J}_3)$.
\end{itemize}
Hence, the eigenvectors of $\hat{\mathbf{F}}_m^{(l)}$ and $\mathbf{\hat{C}}'\mathbf{Q}^{(l)}\mathbf{\hat{B}}$ are the same and then
\[\begin{aligned}
Sp(\hat{\mathbf{F}}_m^{(l)})\ &&=&\ \{m\}\ \cup\ \left\{(1+m)-\lambda,\lambda\in Sp(W^{(l)})\setminus(\{1\}\cup\mathcal{A}_{OUT})\right\}\\
&&\cup&\ \left\{1-\lambda,\lambda\in Sp(\mathbf{Q}^{(l)})\setminus (Sp(W^{(l)})\cup\mathcal{A}_{OUT})\right\},
\end{aligned}\]
which implies $\{\Re e(Sp(\mathcal{D}f(\theta^*)))>1/2\}\equiv\{\Re e(\lambda^{*l})<1/2\}$.

Then, by using analogous arguments of the proof of Theorem~\ref{thm:second_order_leading} for the leading systems,
it can be easily shown that
$$\mathbf{E}[\mathbf{\hat{C}}'\Delta\mathbf{M}^{(l)}_{n}(\Delta\mathbf{M}^{(l)}_{n})^{'}\mathbf{\hat{C}}|\mathcal{F}_{n-1}]\stackrel{a.s.}{\longrightarrow}
\mathbf{\hat{G}}^{(l)},\ \quad\
\mathbf{E}[\mathbf{\hat{C}}'\Delta\mathbf{M}^{(l_1)}_{n}(\Delta\mathbf{M}^{(l_2)}_{n})^{'}\mathbf{\hat{C}}]=0\ \ \ \forall l_1\neq l_2,$$
and for any $\epsilon>0$
\begin{equation*}
\mathbf{E}\left[n\|\mathbf{\hat{C}}'\mathbf{R}^{(l)}_{n}\|^2\mathds{1}_{\{\left\|\mathbf{\hat{Z}}^{(l)}_n-\mathbf{\hat{Z}}^{(l)}_{\infty}\right\|\leq \epsilon\}}\right]\longrightarrow 0.
\end{equation*}
We can then apply Theorem~\ref{ThmCLT} to obtain the CLT with asymptotic variance
$$\hat{\Sigma}^{(l)}\ :=\ \lim_{m\rightarrow\infty}\int_0^{\infty}e^{u(\frac{\mathbf{I}}{2}-\hat{\mathbf{F}}_m^{(l)})}\mathbf{\hat{G}}^{(l)}
e^{u(\frac{\mathbf{I}}{2}-\hat{\mathbf{F}}_m^{(l)})'}du.$$
This concludes the proof.
\endproof


\section*{Acknowledgement}
We thank the anonymous referees for their careful reading of our paper and their
insightful comments and suggestions.

\appendix

\begin{center}
\huge{{\bf Appendix}}
\end{center}

\section{Basic tools of Stochastic Approximation}

We report the recursive procedure defined in~\eqref{SAP} on a filtered probability space
$(\Omega,{\mathcal A},(\mathcal{F}_n)_{n\geq0},\mathbf{P})$, namely
\begin{equation*}\tag{\ref{SAP}}
\forall\, n\geq1,\quad \theta_{n}=\theta_{n-1}-\frac{1}{n}f(\theta_{n-1})+\frac{1}{n}\left(\Delta M_{n}+R_{n}\right),
\end{equation*}
where $f:\mathbb{R}^d\rightarrow \mathbb{R}^d$ is a locally Lipschitz continuous function, $\theta_{n}$ an $\mathcal{F}_{n}$-measurable finite random vector and, for every $n\ge 1$,  $\Delta M_{n}$ is an $\mathcal{F}_{n-1}$-martingale increment and $R_{n}$ is an $\mathcal{F}_n$-adapted remainder term.
\begin{theo}[$A.s.$ convergence with $ODE$ method, see $e.g.$~\cite{Ben, BMP, Duf2, ForPag, KusYin}] \label{ThmODE}
Assume that $f$ is locally Lipschitz, that $$R_n\stackrel{a.s.}{\longrightarrow}0 \quad \mbox{and} \quad
\sup_{n\geq 1}\mathbf{E}\left[\left\|\Delta M_{n}\right\|^2\,|\,\mathcal{F}_{n-1}\right]<+\infty \quad a.s.$$
Then, the set $\Theta^{\infty}$ of its limiting values as $n\rightarrow+\infty$ is $a.s.$ a compact connected set, stable by the flow of
$$ODE_f\equiv\dot{\theta}=-f(\theta).$$
Furthermore, if $\theta^*\in\Theta^{\infty}$ is a uniformly stable equilibrium on $\Theta^{\infty}$ of $ODE_f$, then
$$\theta_n\stackrel{a.s.}{\longrightarrow}\theta^*.$$
\end{theo}

\paragraph{{\sc Comments}} By uniformly stable we mean that
$$\sup_{\theta\in\Theta^{\infty}}\left|\theta(\theta_0,t)-\theta^*\right| \longrightarrow 0\quad\mbox{as} \quad t\rightarrow+\infty,$$
where $\theta(\theta_0,t)_{\theta_0\in \Theta^{\infty}\!\!,t\in \mathbb{R}_+}$ is the flow of $ODE_f$ on $\Theta^{\infty}$.

\medskip
We say that the  function $f$ is $\epsilon${\em-differentiable}, $\epsilon>0$, at $\theta^*$  if
$$f(\theta)=f(\theta^*)+\mathcal{D}f(\theta^*)(\theta-\theta^*)+o(\left\|\theta-\theta^*\right\|^{1+\epsilon}) \quad\mbox{as}\quad \theta\to\theta^*.$$

\begin{theo}[{Rate of convergence see \cite[Theorem 3.III.14 p.131]{Duf2}, for $CLT$ see also $e.g.$~\cite{BMP,KusYin}}]  \label{ThmCLT}
Let $\theta^*$ be an equilibrium point of $\{f = 0\}$. Assume that the function $f$ is differentiable at $\theta^*$ and all the eigenvalues of $\mathcal{D}f(\theta^*)$ have positive real parts. Assume that for some $\delta>0$,
\begin{equation}\label{HypDM1}
	\sup_{n\geq 1}\mathbf{E}\left[\left\|\Delta M_{n}\right\|^{2+\delta}\,|\,\mathcal{F}_{n-1}\right]<+\infty \, a.s.,
 \end{equation}
and \begin{equation}\label{HypDM2}
\mathbf{E}\left[\Delta M_{n}\Delta M_{n}'\,|\,\mathcal{F}_{n-1}\right]\overset{a.s.}{\underset{n\rightarrow+\infty}{\longrightarrow}}\Gamma,
\end{equation}
where $\Gamma\!\in {\mathcal S}^+(d,\mathbb{R})$ (deterministic symmetric  positive matrix) and for an $\epsilon>0$,
\begin{equation}\label{HypReste}
n\mathbf{E}\left[\left\|R_{n}\right\|^2\mathds{1}_{\{\left\|\theta_{n-1}-\theta^*\right\|\leq\epsilon\}}\right]
\underset{n\rightarrow+\infty}{\longrightarrow}0.
\end{equation}
$(a)$ If $\Re e(\lambda_{{\rm min}})>\frac{1}{2}$, where $\lambda_{{\rm min}}$ denotes the eigenvalue of
$\mathcal{D}f(\theta^*)$ with lowest real part, the above a.s. convergence is ruled on the set
$\mathcal{D}f\{\theta_n\rightarrow\theta^*\}$ by the following Central Limit Theorem
$$
\sqrt{n}\left(\theta_n-\theta^*\right)\overset{{\mathcal L}}{\underset{n\rightarrow\infty}{\longrightarrow}}{\mathcal N}
\left(0,\Sigma\right) \quad
\mbox{with} \quad \Sigma:=\displaystyle\int_0^{+\infty}e^{\left(I_d/2-\mathcal{D}f(\theta^*)\right)u}
\Gamma e^{\left(I_d/2-\mathcal{D}f(\theta^*)\right)^{'}u}du.$$
\noindent $(b)$ If $\Re e(\lambda_{{\rm min}})=\frac{1}{2}$, then
$$\sqrt{\frac{n}{\log n}}\left(\theta_n-\theta^*\right)\overset{{\mathcal L}}
{\underset{n\rightarrow\infty}{\longrightarrow}}{\mathcal N}(0,\Sigma)\quad\mbox{as $n\to+\infty$}.$$
$(c)$  If $ \Re e(\lambda_{{\rm min}})\!\in(0,\frac{1}{2})$, then $n^{\Re e(\lambda_{{\rm min}})}\left(\theta_n-\theta^*\right)$ $a.s.$ converges as $n\to+\infty$ towards a finite random variable.
\end{theo}


\begin{thebibliography}{10}

\bibitem{AleMaySec1}
G.~Aletti, C.~May, and P.~Secchi.
\newblock On the distribution of the limit proportion for a two-color, randomly
  reinforced urn with equal reinforcement distributions.
\newblock {\em Adv. in Appl. Probab.}, 39(3):690--707, 2007.

\bibitem{AleMaySec2}
G.~Aletti, C.~May, and P.~Secchi.
\newblock A central limit theorem, and related results, for a two-color
  randomly reinforced urn.
\newblock {\em Adv. in Appl. Probab.}, 41(3):829--844, 2009.

\bibitem{Ath}
K.~B. Athreya.
\newblock On a characteristic property of {P}\'{o}lya's urn.
\newblock {\em Studia Sci. Math. Hungar.}, 4:31--35, 1969.

\bibitem{AthKar2}
K.~B. Athreya and S.~Karlin.
\newblock Embedding of urn schemes into continuous time {M}arkov branching
  processes and related limit theorems.
\newblock {\em Ann. Math. Statist.}, 39:1801--1817, 1968.

\bibitem{BaiHu}
Z.~D. Bai and F.~Hu.
\newblock Asymptotic theorems for urn models with nonhomogeneous generating
  matrices.
\newblock {\em Stoch. Proc. Appl.}, 80(1):87--101, 1999.

\bibitem{BaiHu2}
Z.-D. Bai and F.~Hu.
\newblock Asymptotics in randomized urn models.
\newblock {\em Ann. Appl. Probab.}, 15(1B):914--940, 2005.

\bibitem{Ben}
M.~Bena{\"{\i}}m.
\newblock Dynamics of stochastic approximation algorithms.
\newblock In {\em S\'eminaire de {P}robabilit\'es, {XXXIII}}, volume 1709 of
  {\em Lecture Notes in Math.}, pages 1--68. Springer, Berlin, 1999.

\bibitem{BenBenCheLim}
M.~Bena{\"{\i}}m, I.~Benjamini, J.~Chen, and Y.~Lima.
\newblock A generalized {P}\'{o}lya's urn with graph based interactions.
\newblock {\em Random Structures Algorithms}, 46(4):614--634, 2015.

\bibitem{BMP}
A.~Benveniste, M.~M{\'e}tivier, and P.~Priouret.
\newblock {\em Adaptive algorithms and stochastic approximations}, volume~22 of
  {\em Applications of Mathematics (New York)}.
\newblock Springer-Verlag, Berlin, 1990.
\newblock Translated from the French by Stephen S. Wilson.

\bibitem{CheLuc}
J.~Chen and C.~Lucas.
\newblock A generalized {P}\'{o}lya's urn with graph based interactions:
  convergence at linearity.
\newblock {\em Electron. Commun. Probab.}, 19:no. 67, 13, 2014.

\bibitem{CirGalHus}
P.~Cirillo, M.~Gallegati, and J.~H{\"u}sler.
\newblock A {P}\'olya lattice model to study leverage dynamics and contagious
  financial fragility.
\newblock {\em Adv. Complex Syst.}, 15(suppl. 2):1250069, 26, 2012.

\bibitem{CriDaiPraMin}
I.~Crimaldi, P.~Dai~Pra, and I.~G. Minelli.
\newblock Fluctuation theorems for synchronization of interacting {P}\'olya's
  urns.
\newblock {\em Stochastic Process. Appl.}, 126(3):930--947, 2016.

\bibitem{DaiPraLouisMin}
P.~Dai~Pra, P.-Y. Louis, and I.~G. Minelli.
\newblock Synchronization via interacting reinforcement.
\newblock {\em J. Appl. Probab.}, 51(2):556--568, 2014.

\bibitem{Duf2}
M.~Duflo.
\newblock {\em Random iterative models}, volume~34 of {\em Applications of
  Mathematics (New York)}.
\newblock Springer-Verlag, Berlin, 1997.
\newblock Translated from the 1990 French original by Stephen S. Wilson and
  revised by the author.

\bibitem{DurFloLi}
S.~D. Durham, N.~Flournoy, and W.~Li.
\newblock A sequential design for maximizing the probability of a favourable
  response.
\newblock {\em Canad. J. Statist.}, 26(3):479--495, 1998.

\bibitem{EffPol}
F.~Eggenberger and G.~{P}\'{o}lya.
\newblock {\"U}ber die statistik verketteter vorg{\"a}nge.
\newblock {\em Math. Mech.}, 3:279--289, 1923.

\bibitem{ForPag}
J.-C. Fort and G.~Pag{\`e}s.
\newblock Convergence of stochastic algorithms: from the {K}ushner-{C}lark
  theorem to the {L}yapounov functional method.
\newblock {\em Adv. in Appl. Probab.}, 28(4):1072--1094, 1996.

\bibitem{Fri}
B.~Friedman.
\newblock A simple urn model.
\newblock {\em Comm. Pure Appl. Math.}, 2:59--70, 1949.

\bibitem{HilLanSud1}
B.~M. Hill, D.~Lane, and W.~Sudderth.
\newblock A strong law for some generalized urn processes.
\newblock {\em Ann. Probab.}, 8(2):214--226, 1980.

\bibitem{HilLanSud2}
B.~M. Hill, D.~Lane, and W.~Sudderth.
\newblock Exchangeable urn processes.
\newblock {\em Ann. Probab.}, 15(4):1586--1592, 1987.

\bibitem{JohKot}
N.~L. Johnson and S.~Kotz.
\newblock {\em Urn models and their application}.
\newblock John Wiley \& Sons, New York-London-Sydney, 1977.
\newblock An approach to modern discrete probability theory, Wiley Series in
  Probability and Mathematical Statistics.

\bibitem{KusYin}
H.~J. Kushner and G.~G. Yin.
\newblock {\em Stochastic approximation and recursive algorithms and
  applications}, volume~35 of {\em Applications of Mathematics (New York)}.
\newblock Springer-Verlag, New York, second edition, 2003.
\newblock Stochastic Modelling and Applied Probability.

\bibitem{LarPag}
S.~Laruelle and G.~Pag{\`e}s.
\newblock Randomized urn models revisited using stochastic approximation.
\newblock {\em Ann. Appl. Probab.}, 23(4):1409--1436, 2013.

\bibitem{Lau}
M.~Launay.
\newblock Interacting urn models.
\newblock {\em arXiv.org}, math.PR, 2011.

\bibitem{LauLim}
M.~Launay and V.~Limic.
\newblock Generalized interacting urn models.
\newblock {\em arXiv.org}, math.PR, 2012.

\bibitem{MarVal}
M.~Marsili and A.~Valleriani.
\newblock {Self organization of interacting {P}\'{o}lya urns}.
\newblock {\em {European Physical Journal B}}, {3}({4}):{417--420}, {JUN}
  {1998}.

\bibitem{NalParTos}
G.~Naldi, L.~Pareschi, and G.~Toscani, editors.
\newblock {\em Mathematical modeling of collective behavior in socio-economic
  and life sciences}.
\newblock Modeling and Simulation in Science, Engineering and Technology.
  Birkh\"auser Boston, Inc., Boston, MA, 2010.

\bibitem{Nor}
J.~R. Norris.
\newblock {\em Markov chains}, volume~2 of {\em Cambridge Series in Statistical
  and Probabilistic Mathematics}.
\newblock Cambridge University Press, Cambridge, 1998.
\newblock Reprint of 1997 original.

\bibitem{PagSec}
A.~M. Paganoni and P.~Secchi.
\newblock Interacting reinforced-urn systems.
\newblock {\em Adv. in Appl. Probab.}, 36(3):791--804, 2004.

\bibitem{Pem}
R.~Pemantle.
\newblock A survey of random processes with reinforcement.
\newblock {\em Probab. Surv.}, 4:1--79, 2007.

\bibitem{Pri}
A.~Prindle, P.~Samayoa, I.~Razinkov, T.~Danino, L.~S. Tsimring, and J.~Hasty.
\newblock A sensing array of radically coupled genetic biopixels.
\newblock {\em Nature}, 481:39--44, 2011.

\bibitem{Shi}
A.~N. Shiryaev.
\newblock {\em Probability}, volume~95 of {\em Graduate Texts in Mathematics}.
\newblock Springer-Verlag, New York, second edition, 1996.
\newblock Translated from the first (1980) Russian edition by R. P. Boas.

\bibitem{Smy}
R.~T. Smythe.
\newblock Central limit theorems for urn models.
\newblock {\em Stoch. Proc. Appl}, 65(1):115--137, 1996.

\bibitem{Will}
D.~Williams.
\newblock {\em Probability with martingales}.
\newblock Cambridge Mathematical Textbooks. Cambridge University Press,
  Cambridge, 1991.

\bibitem{YakTso}
A.~Y. Yakovlev and A.~D. Tsodikov.
\newblock {\em Stochastic Models of Tumor Latency and their Biostatistical
  Applications}.
\newblock World Scientific, Singapore, 1996.

\end{thebibliography}
\end{document}